\documentclass[a4paper]{article}

\usepackage[utf8]{inputenc}
\usepackage[UKenglish]{babel}
\usepackage[T1]{fontenc}           	
\usepackage{amsmath}
\usepackage{amssymb}
\usepackage{amsfonts}          		
\usepackage[all]{xy}
\usepackage{amscd}
\usepackage{amsthm}
\usepackage{graphics}
\usepackage{color}
\usepackage{mathtools}
\usepackage[a4paper]{geometry}
\usepackage{bbm}
\usepackage{mathrsfs}
\usepackage{lipsum}
\usepackage{indentfirst} 
\usepackage{latexsym} 
\usepackage{xcolor}
\usepackage{cite}
\usepackage[hidelinks]{hyperref}
\usepackage{algorithm}
\usepackage{algpseudocode}
\usepackage{enumitem}
\usepackage{multirow}
\usepackage{multicol}
\usepackage{xspace}

\newtheorem{thm}{Theorem}[section]

\newtheorem{lem}[thm]{Lemma}

\newtheorem{rem}{Remark}

\newtheorem{thm*}{Theorem}

\newtheorem{lem*}{Lemma}

\newcommand{\R}{\mathbb{R}}

\newcommand{\io}{\mathcal{O}}

\newcommand{\symm}{\textnormal{Sym}}
\newcommand{\skeww}{\textnormal{Skew}}
\newcommand{\rank}{\textnormal{rank}}
\newcommand{\trace}{\textnormal{trace}}
\newcommand{\Matlab}{\textsc{Matlab}\xspace}

\DeclareMathOperator*{\argmin}{arg\,min}

\numberwithin{figure}{section}
\numberwithin{table}{section}

\definecolor{brightpink}{rgb}{1.0, 0.0, 0.5}

\newcommand{\revise}[1]{{\color{black} #1}}

\newcommand{\reviseLAA}[1]{{\color{black} #1}}

\title{Minimum-norm solutions of the non-symmetric semidefinite Procrustes problem} 
\author{
Nicolas Gillis\thanks{Department of Mathematics and Operational Research, University of Mons, Mons, Belgium. We acknowledge  the support by the European Union (ERC consolidator, eLinoR, no 101085607). Email: Nicolas.GILLIS@umons.ac.be, Stefano.SICILIA@umons.ac.be. } 
\and Stefano Sicilia\footnotemark[1] $^,$\thanks{Part of this work was done when SS was with the Division of Mathematics, Gran Sasso Science Institute, L'Aquila, Italy. SS  is member of the Gruppo Nazionale Calcolo Scientifico-Istituto Nazionale di Alta Matematica (GNCS-INdAM).}}
\date{}

\begin{document}
 \maketitle
 
 \begin{abstract}
  Given two matrices $X,B\in \R^{n\times m}$ and a set $\mathcal{A}\subseteq \R^{n\times n}$, a Procrustes problem  consists in finding a matrix $A \in \mathcal{A}$ such that the Frobenius norm of $AX-B$ is minimized. When $\mathcal{A}$ is the set of the matrices whose symmetric part is positive semidefinite, we obtain the so-called non-symmetric positive semidefinite Procrustes (NSPSDP) problem.  The NSPSDP problem arises in the estimation of compliance or stiffness matrix in solid and elastic structures. 
  If $X$ has rank $r$, Baghel et al.\ (Lin.\ Alg.\ Appl., 2022) proposed a three-step semi-analytical approach: 
  (1) construct a reduced NSPSDP problem in dimension $r\times r$, 
  (2) solve the reduced problem by means of a fast gradient method with a linear rate of convergence, and 
  (3) post-process the solution of the reduced problem to construct a solution of the larger original NSPSDP problem. In this paper, we revisit this approach of  Baghel et al.\ and identify an unnecessary assumption used by the authors leading to cases where their algorithm cannot attain a minimum and produces solutions with unbounded norm.  
  In fact, revising the post-processing phase of their semi-analytical approach, we show that the infimum of the NSPSDP problem is always attained, and we show how to compute a minimum-norm solution. 
  We also prove that the symmetric part of the computed solution has minimum rank bounded by $r$, and that the skew-symmetric part has rank bounded by $2r$. Several numerical examples show the efficiency of this algorithm, both in terms of computational speed and of finding optimal minimum-norm solutions. 
 \end{abstract}
 
 \vspace{0.1cm}
 

 \section{Introduction}
 
 A structured mapping problem associated to the set $\mathcal{A}\subseteq \R^{n\times n}$ consists in finding a matrix $A\in \mathcal{A}$ such that $AX=B$, where $X,B\in \R^{n\times m}$ are two fixed matrices; see, e.g.,  \cite{adhikari2008backward, mackey2008structured} and the references therein. 
 The corresponding Procrustes problem aims to solve a mapping structured problem when it is inconsistent, that is,  when there is no matrix $A\in \mathcal{A}$ such that $AX=B$, leading to the following optimization problem 
 \begin{equation}
  \label{prob:proc}
  \inf_{A\in \mathcal{A}} \|AX-B\|, \quad \text{ where $\| \cdot \|$ is some norm. }
 \end{equation} 
 Many constraints have been studied for the set $\mathcal{A}$, such as the orthogonal matrices arising in computer graphics \cite{green1952orthogonal, gower2004procrustes}, Jordan algebra and Lie algebra for which Adhikari and Alam have provided an analytical solution \cite{adhikari2016structured},  positive semidefiniteness in \cite{brock1968optimal, suffridge1993approximation} arising in elastic structures, and rank-$r$ matrices with additional constraints in~\cite{li2023generalized},  
 to name just a few.

 In this paper we focus on a set $\mathcal{A}$ requiring a semidefinite property of the matrix. There are two main cases and the set $\mathcal{A}$ is one of the following 
 \[
  \mathcal{S}_\succeq^n:=\{M\in \R^{n\times n} : M\succeq 0\}, \quad  \text{ or } \quad  
  \mathcal{N}_\succeq^n:=\{M\in \R^{n\times n} : (M+M^\top)\succeq 0\},
 \]
 where $M\succeq N$ means that $M-N$ is positive semidefinite, and $0$ is the zero matrix of appropriate dimension. 
 The set $\mathcal{S}_\succeq^n$ denotes the set of all symmetric and positive semidefinite matrices, while the set $\mathcal{N}_\succeq^n$ contains all the matrices whose symmetric part is positive semidefinite. Given $X,B\in \R^{n\times m}$, the positive semidefinite Procrustes (PSDP) problem is 
 \begin{equation}
  \tag{\ensuremath{\mathcal{P}}}
  \label{prob:PSDP} 
  \inf_{A \in \mathcal{S}_\succeq^n} \|A X-B\|^2_F,
 \end{equation}
 while the non-symmetric positive semidefinite Procrustes (NSPSDP) problem is
 \begin{equation}
 \tag{\ensuremath{\mathcal{NP}}}
  \label{prob:NSPSDP}
  \inf_{A \in \mathcal{N}_\succeq^n} \|A X-B\|^2_F,
 \end{equation}
 where $\|\cdot\|_F$ denotes the Frobenius norm. These problems are both convex, since they minimize a convex function over a convex set. 
 While \eqref{prob:PSDP} is well-known in the literature and  has been studied in many works \cite{brock1968optimal,woodgate1996least,krislock2003numerical,gillis2018semi,jingjing2019solution},  \eqref{prob:NSPSDP} is more recent and has been considered in \cite{krislock2003numerical,krislock2004local,baghel2022non}. 
 Both PSDP and NSPSDP problems arise in the estimate of local compliance matrices during deformable object in various engineering applications; see \cite{krislock2004local,baghel2022non}. The compliance matrix represents the Green function matrix that corresponds to the unknown $A$ in \eqref{prob:NSPSDP}. The estimation of $A$ is based on the observations of the forces, stored in the columns of $X$, and their associated displacements, stored in the columns of $B$, for a total of $m$ measurements. When the deformable object is a passive object, that is, it does not generate energy in deformation, the physical properties of the model impose that $p^\top A p\geq 0$, where $p$ is the point load, which is implied by $A \in \mathcal{N}^n_\succeq$. Thus, when estimating $A$ under noisy measurements, it is crucial to impose the constraint $A \in \mathcal{N}^n_\succeq$ in the optimization problem; \reviseLAA{ see \cite{lang2002acquisition} and \cite[Section III]{krislock2004local} for more details about the physical model}.

 \subsection{Contributions and outline of the paper}

For both \eqref{prob:PSDP} and \eqref{prob:NSPSDP}, the original $n\times m$ dimensional problem can be reduced to a  problem of dimension $r\times r$, where $r$ is the rank of $X$; see  
  \cite{gillis2018semi} for  \eqref{prob:PSDP},  and \cite{baghel2022non} for \eqref{prob:NSPSDP}.  This is referred to as a semi-analytical approach because the smaller $r \times r$ dimensional problem still needs to be solved with iterative methods, and \cite{gillis2018semi, baghel2022non} relied on an optimal first-order method with linear convergence. 
It is important to note that the solution of~\eqref{prob:PSDP} is not always attained~\cite{gillis2018semi}; see Section~\ref{sec:PSDP} where we recall this result.  
The approach in~\cite{baghel2022non} follows closely that of \cite{gillis2018semi}, and concludes that the solution of~\eqref{prob:NSPSDP} is also not always attained. However, it turns out this result is not correct, unless one uses an additional \revise{constraint} within the NSPSDP problem, which is not explicitly stated in~\cite{baghel2022non}.  
 
 In this work, inspired by the approach of  \cite{baghel2022non}, we provide an updated theorem about the NSPSDP problem, namely, we show that the infimum is always attained and provide a semi-analytical solution  (Theorem~\ref{thm:solNSPSDP}).  
 Moreover, we show how to compute a minimum-norm solution for which we also prove low-rank properties, namely the symmetric part has rank at most $r$, and the skew-symmetric part at most $2r$ (Theorem~\ref{thm:minnormNSPSDP}). 
 This allows us to propose a new algorithmic approach for the NSPSDP problem, which always provides optimal and minimum-norm solutions, as opposed to~\cite{baghel2022non}. 
  
 The paper is organized as follows. 
 In Section~\ref{sec:overview}, we recall the semi-analytical methods proposed in \cite{gillis2018semi} for solving \eqref{prob:PSDP} and in \cite{baghel2022non} for solving \eqref{prob:NSPSDP}. We discuss in detail the main issue of the latter approach, namely the introduction of an unnecessary additional constraint in the problem, which leads to non-optimal and unbounded solutions. 
 In Section~\ref{sec:newappr}, we show how to fix this issue, and prove that the infimum of the NSPSDP problem is always attained. \revise{We also discuss} how to compute the minimum-norm solution and show its low-rank properties. 
 In Section~\ref{sec:algor}, we develop an algorithm relying on the theoretical results from Section~\ref{sec:newappr}. 
 Section~\ref{sec:numexp} shows the effectiveness of this new algorithm compared to that of~\cite{baghel2022non} in several numerical examples.

 \section{An overview of  Procrustes problems with positive semidefinite constraints} 
 \label{sec:overview}

 In this section, we briefly recall the existing approaches by Gillis and Sharma~\cite{gillis2018semi} for solving \eqref{prob:PSDP}, and by Baghel et al.\ \cite{baghel2022non} for solving \eqref{prob:NSPSDP}. 
For the latter, we provide a revised version of their main theorem, stating explicitely the hidden constraint in their result.

 \subsection{The PSDP problem} \label{sec:PSDP}

 Let us briefly recall the existence of optimal solutions in the PSDP problem; this will allow us to compare and shed light on the NSPSD case. 

 Before doing this, let us recall an important result about \revise{positive semidefinite}    characterization that will be used throughout the paper.  
 \begin{lem}\cite{albert1969conditions} 
  \label{lem:blockPD}
  Let $R\in \R^{n\times n}$ be a symmetric matrix partitioned as
  \[
   R=\begin{pmatrix}
      B & C^\top \\
      C & D \\
     \end{pmatrix}.
  \]
  Then $R\succeq 0$ if and only if the following three conditions hold: (1) $B\succeq 0$,  
(2) $\ker(B)\subseteq \ker(C)$, and 
  (3) $D-CB^\dagger C^\top\succeq 0$, where $B^\dagger$ denotes the Moore-Penrose pseudo-inverse of $B$. 
 \end{lem} 
 
 When the matrix $X$ has full row rank, the solution to problem \eqref{prob:PSDP} is unique, as ensured by the next lemma.
 \begin{lem}\cite[Theorem 2.2]{woodgate1996least}
  \label{lem:uniquePSDP}
  Let $X,B\in \R^{n \times m}$ and assume that $X$ has rank $n$. Then the infimum of problem \eqref{prob:PSDP} is attained for a unique solution $A_\star \in \mathcal{S}_\succeq^n$.
  \begin{proof}
  This follows from the fact that the problem is strongly convex. 
  \end{proof}
 \end{lem} 
 Now we state the major result of \cite{gillis2018semi}, where a minimum-norm solution of the Procrustes problem \eqref{prob:PSDP} is provided via a semi-analytical approach. 
 \begin{thm}
  \label{thm:solPSDP}
  Let $X,B\in \R^{n\times m}$ and consider the SVD 
  \[
   X=U\Sigma V^\top:=
   \begin{pmatrix}
    U_1 & U_2 \\
   \end{pmatrix}
   \begin{pmatrix}
    \Sigma_1 & 0 \\
    0 & 0 \\
   \end{pmatrix}
   \begin{pmatrix}
    V_1 & V_2 \\
   \end{pmatrix}^\top,
  \]
  where $\Sigma_1\in \R^{r\times r}$  and $r = \rank(X)$. Then
  \begin{equation} \label{eq:infPSDP}
   \inf_{A \in \mathcal{S}_\succeq^n} \|A X-B\|^2_F=\min_{A_{11}\in \mathcal{S}_\succeq^r}\|A_{11}\Sigma_1-U_1^\top B V_1\|_F^2+\|BV_2\|_F^2.
  \end{equation}
  Moreover, letting 
  \[
   \hat{A}_{11}:=\argmin_{A_{11}\in \mathcal{S}_\succeq^r}\|A_{11}\Sigma_1-U_1^\top B V_1\|_F^2,
  \]
  and $Z:=U_2^\top B V_1 \Sigma_1^{-1}$, it holds that 
  \begin{enumerate}
   \item If $\ker(\hat{A}_{11})\subseteq \ker(Z)$, the matrix 
   \[
    A_\star=U_1\hat{A}_{11}U_1^\top+U_2ZU_1^\top+U_1Z^\top U_2^\top+U_2KU_2^\top,
   \]
   where $K=Z\hat{A}_{11}^\dagger Z^\top$, is the unique matrix that attains the infimum in~\eqref{eq:infPSDP} and whose rank, Frobenius norm and spectral norm are minimum. 
   
   \item Otherwise the infimum is not attained, and we can construct a solution  \reviseLAA{$A_\star^{(\varepsilon)}\in \mathcal{S}_\succeq^n$} with minimal rank, minimum Frobenius norm and minimum spectral norm, such that 
   \[
    \|A_\star^{(\varepsilon)} X-B\|^2_F\leq \inf_{A\in \mathcal{S}_\succeq^n}\|AX-B\|^2_F+\varepsilon, 
   \]
   for any $\varepsilon > 0$ sufficiently small. 
  \end{enumerate}
  \begin{proof}
   We refer to \cite[Theorem 1]{gillis2018semi} and 
   \cite[Corollary 1]{gillis2018semi} for the proof.  
  \end{proof}
 \end{thm}
 
 \reviseLAA{
 \begin{rem}
  \label{rem:epslong}
  The explicit expressions for $\varepsilon$ and $A_\star^{(\varepsilon)}$ in Theorem \ref{thm:solPSDP} in the case the infimum is not attained can be found in \cite{gillis2018semi}. These formulas are analogous to the ones of Theorem~\ref{thm:preNSPSDP} presented 
  in the next section.          
 \end{rem}
 }

 \subsection{The NSPSDP problem} \label{sec:NSPSDP} 
 
 The result of Theorem~\ref{thm:solPSDP} for problem \eqref{prob:PSDP} has been adapted by Baghel et al.~\cite{baghel2022non} for~\eqref{prob:NSPSDP}, but the approach they proposed contains some incorrect statements. Before showing the issues in the main result of \cite{baghel2022non}, we start by recalling an analogous result to Lemma \ref{lem:uniquePSDP} for the NSPSDP problem, which guarantees a unique solution when $X$ has full row rank. 
 \begin{lem} \cite[Theorem 2.4.6]{krislock2003numerical} 
  \label{lem:uniqueNSPSDP}
  Let $X,B\in \R^{n \times m}$ and assume that $X$ has rank $n$. Then the infimum of problem \eqref{prob:NSPSDP} is attained for a unique solution $A_\star \in \mathcal{N}_\succeq^n$.
  \begin{proof}
   This follows from the fact that the problem is strongly convex. 
  \end{proof}
 \end{lem}
 
 Now we discuss the major result in \cite{baghel2022non}. 
 Let us define the set
 \begin{equation}
  \label{eq:cset}
  \mathcal{C}=\left\{ \left(
  \begin{matrix}
   C_{11} & 0 \\
   C_{12} & C_{22} \\
  \end{matrix}\right) : C_{11}\in \R^{r\times r}, \ 
  C_{12} \in \R^{(n-r) \times r}, \
  C_{22} \in \R^{(n-r)\times (n-r)}  \right\}\subseteq  \R^{n\times n}, 
 \end{equation}
 where there is a $r\times (n-r)$ block of zeros at the top right of any matrix in $\mathcal{C}$.  
For a given matrix $U\in \R^{n\times n}$, let us also introduce the sets
 \[
   U\mathcal{C}U^\top=\left\{ UCU^\top: C\in \mathcal{C}\right\}, \qquad \text{ and }  \qquad \mathcal{N}_\succeq^{U\mathcal{C}U^\top}=\mathcal{N}_\succeq^n\cap U\mathcal{C}U^\top. 
 \]
 The main issue with \cite[Theorem 2]{baghel2022non} is that it considers a different version of the NSPSDP problem, where it is implicitly assumed that the solution $A$ belongs to $U\mathcal{C}U^\top$, 
 where  the matrix $U$ comes from the SVD of $X=U\Sigma V^\top$. 
 Hence the problem considered  in \cite{baghel2022non} is actually 
 \begin{equation}
  \tag{\ensuremath{\mathcal{NPC}}}
  \label{prob:preNSPSDP}
  \inf_{A \in \mathcal{N}_\succeq^{U\mathcal{C}U^\top}} \|A X-B\|^2_F,
 \end{equation}
 which is different than the original NSPSDP problem \eqref{prob:NSPSDP}. Although the infimums  of~\eqref{prob:preNSPSDP} and \eqref{prob:NSPSDP} coincide, the optimal solution of~\eqref{prob:preNSPSDP} might not be attained in some cases, see Theorem~\ref{thm:preNSPSDP} below, while that of \eqref{prob:NSPSDP} always is; see Theorem~\ref{thm:solNSPSDP} in the next section. 
 
 In the following we provide a revised version of \cite[Theorem 2]{baghel2022non} with this additional constraint which is  not present in the original paper but was  accidentally imposed in the proof of the result. 
 \begin{thm} \cite[Theorem 2, revised]{baghel2022non}
  \label{thm:preNSPSDP}
  Let $X,B\in \R^{n\times m}$ and consider the SVD 
  \begin{equation}   \label{eq:revisedNSPSDP}
   X=U\Sigma V^\top:=
   \begin{pmatrix}
    U_1 & U_2 \\
   \end{pmatrix}
   \begin{pmatrix}
    \Sigma_1 & 0 \\
    0 & 0 \\
   \end{pmatrix}
   \begin{pmatrix}
    V_1 & V_2 \\
   \end{pmatrix}^\top,
  \end{equation}
  where $\Sigma_1\in \R^{r\times r}$ has full rank. Then
  \[
   \inf_{A \in \mathcal{N}_\succeq^{U\mathcal{C}U^\top}} \|A X-B\|^2_F=\inf_{A \in \mathcal{N}_\succeq^n} \|A X-B\|^2_F=\min_{A_{11}\in \mathcal{N}_\succeq^r}\|A_{11}\Sigma_1-U_1^\top B V_1\|_F^2+\|BV_2\|_F^2.
  \]
  Moreover, let us introduce
  \[
   \hat{A}_{11}=\argmin_{A_{11}\in \mathcal{N}_\succeq^r}\|A_{11}\Sigma_1-U_1^\top B V_1\|_F^2:=\hat{H}_{11}+\hat{S}_{11},
  \]
  where $\hat{H}_{11}$ and $\hat{S}_{11}$ are, respectively, the symmetric and skew-symmetric part of $\hat{A}_{11}$. Let $Z=U_2^\top B V_1\Sigma_1^{-1}$, the following holds:
  \begin{enumerate} 
   \item If $\ker(\hat{H}_{11})\subseteq \ker(Z)$, then a solution of the problem~\eqref{prob:preNSPSDP} is
   \[
    A_\star = U
    \begin{pmatrix}
     \hat{A}_{11} & 0 \\
     Z & K+R \\
\end{pmatrix} U^\top=U_1\hat{A}_{11}U_1^\top+U_2ZU_1^\top+U_2(K+R)U_2^\top,
   \]
   where $R\in \R^{(n-r)\times (n-r)}$ is skew-symmetric and $K\in \R^{(n-r)\times (n-r)}$ is symmetric such that
   \[
    K \succeq \frac{1}{4}ZH_{11}^\dagger Z^\top.
   \]
Since the infimum of~\eqref{prob:NSPSDP} and~\eqref{prob:preNSPSDP} \reviseLAA{coincide (see \eqref{eq:revisedNSPSDP})}, $A_\star$ is also an optimal solution of~\eqref{prob:NSPSDP}. 
   
   \item Otherwise the infimum \reviseLAA{of \eqref{prob:preNSPSDP}} is not attained, and we can construct a solution  \reviseLAA{$A_\star^{(\varepsilon)}\in \mathcal{N}_\succeq^n$} such that
   \[
    \|A_\star^{(\varepsilon)} X-B\|^2_F\leq \inf_{A\in \mathcal{N}_\succeq^n}\|AX-B\|^2_F+\varepsilon 
   \]
for any     
   \[
    0<\varepsilon<
    \begin{dcases}
     \min\{1,\|\hat{A}_{11}\Sigma_1-U_1^\top B V_1\|_F^2\} & \textnormal{if } \|\hat{A}_{11}\Sigma_1-U_1^\top B V_1\|_F^2\neq 0, \\
     1 & \textnormal{otherwise.}
    \end{dcases}
   \]
   Consider the two SVDs
   \[
    \hat{H}_{11}=
   \begin{pmatrix}
    W_1 & W_2 \\
   \end{pmatrix}
   \begin{pmatrix}
    \Lambda & 0 \\
    0 & 0 \\
   \end{pmatrix}
   \begin{pmatrix}
    W_1 & W_2 \\
   \end{pmatrix}^\top, \qquad
   \hat{H}^{(\varepsilon)}_{11}=
   \begin{pmatrix}
    W_1 & W_2 \\
   \end{pmatrix}
   \begin{pmatrix}
    \Lambda & 0 \\
    0 & \Gamma^{(\varepsilon)} \\
   \end{pmatrix}
   \begin{pmatrix}
    W_1 & W_2 \\
   \end{pmatrix}^\top,
   \] 
   where $\Lambda\in \R^{s\times s}$ is invertible and $\Gamma^{(\varepsilon)}\in \R^{r-s\times r-s}$ is the diagonal matrix with entries $\varepsilon/\beta$, with
   \[
    \beta=
    \begin{dcases}
     4\sqrt{r-s}\|\Sigma_1\|_F\|\hat{A}_{11}\Sigma_1-U_1^\top B V_1\|_F & \textnormal{if } \|\hat{A}_{11}\Sigma_1-U_1^\top B V_1\|_F^2\neq 0, \\
     4\sqrt{r-s}\|\Sigma_1\|_F & \textnormal{otherwise.}
    \end{dcases}
   \]
   Defining  $\hat{A}^{(\varepsilon)}_{11}=\hat{H}^{(\varepsilon)}_{11}+\hat{S}_{11}$, we construct $A_\star^{(\varepsilon)}$ as follows 
   \begin{equation} \label{sol:unattained}
    A_\star^{(\varepsilon)}=U
    \left(\begin{matrix}
     \hat{A}_{11} & 0 \\
     Z & K^{(\varepsilon)}+R \\  \end{matrix}\right)U^\top=U_1\hat{A}_{11}U_1^\top+U_2ZU_1^\top+U_2(K^{(\varepsilon)}+G)U_2^\top,
   \end{equation} 
   where $R\in \R^{(n-r)\times (n-r)}$ is skew-symmetric and $K^{(\varepsilon)} \succeq \frac{1}{4}Z(\hat{H}^{(\varepsilon)}_{11})^{-1} Z^\top$. 
  
  \end{enumerate}
  \begin{proof}
   We report the proof of \cite[Theorem 2]{baghel2022non} until the additional constraint $A\in U\mathcal{C}U^\top$ arises, where the set $\mathcal{C}$ is defined in \eqref{eq:cset}. Let $A\in \mathcal{N}_\succeq^n$ and set
   \[
    \hat{A}:=U^\top AU=U^\top
    \begin{pmatrix}
     H_{11} & H_{21}^\top \\
     H_{21} & H_{22} \\
    \end{pmatrix}U+U^\top
    \begin{pmatrix}
     S_{11} & -S_{21}^\top \\
     S_{21} & S_{22} \\
    \end{pmatrix}U:=H+S,
   \]
   where $H$ and $S$ are, respectively, the symmetric and skew-symmetric part of $\hat{A}$ that have been decomposed and partitioned in blocks such that $H_{11},S_{11}\in \R^{r\times r}$. By definition $S_{11}$ and $S_{22}$ are skew-symmetric and, since $H\succeq 0$, Lemma \ref{lem:blockPD} implies that
   \[
    H_{11}\succeq 0, \quad \ker(H_{11})\subseteq\ker(H_{21}), \quad H_{22}-H_{21}H_{11}^\dagger H_{21}^\top\succeq 0.
   \]
   We rewrite the objective function and we get
   \begin{align*}
    & \|AX-B\|_F^2 = \|U\hat{A}U^\top X-B\|_F^2=\|(H+S)U^\top X-U^\top B\|_F^2 
   \\
    & \quad = \|(H_{11}+S_{11})U_1^\top X-U_1^\top B\|_F^2+\|(H_{21}+S_{21})U_1^\top X-U_2^\top B\|_F^2 
   \\
     & \quad = \left\|\left(\begin{matrix}
      (H_{11}+S_{11})\Sigma_1-U_1^\top BV_1 & -U_1^\top BV_2
     \end{matrix}\right)\right\|_F^2+
     \left\|\left(\begin{matrix}
      (H_{21}+S_{21})\Sigma_1-U_2^\top BV_1 & -U_2^\top BV_2
     \end{matrix}\right)\right\|_F^2 
   \\
     & \quad = \|(H_{11}+S_{11})\Sigma_1-U_1^\top BV_1\|_F^2+\|(H_{21}+S_{21})\Sigma_1-U_1^\top BV_1\|_F^2+\|U_1^\top BV_2\|_F^2+\|U_2^\top BV_2\|_F^2,
   \end{align*} 
   which implies
   \begin{equation}
    \label{eq:minpreNSPSDP}
    \|AX-B\|_F^2=\|(H_{11}+S_{11})\Sigma_1-U_1^\top BV_1\|_F^2+\|(H_{21}+S_{21})\Sigma_1-U_2^\top BV_1\|_F^2+\|BV_2\|_F^2.
   \end{equation}
   By taking the infimum of \eqref{eq:minpreNSPSDP}, $\inf_{A\in \mathcal{N}_\succeq^n}\|AX-B\|_F^2$ is equal to 
   \[
    \inf_{\substack{H_{11}\succeq 0,S_{11}=-S_{11}^\top, \\ \ker(H_{11})\subseteq \ker(H_{21})}} \|(H_{11}+S_{11})\Sigma_1-U_1^\top BV_1\|_F^2+\|(H_{21}+S_{21})\Sigma_1-U_2^\top BV_1\|_F^2+\|BV_2\|_F^2. 
   \]
   We note that in this last infimum the matrices $H_{22}$ and $S_{22}$ do not appear, meaning that the conditions $S_{22}=-S_{22}^\top$ and $H_{22}-H_{21}H_{11}^\dagger H_{21}^\top\succeq 0$ can be ignored. By dropping the condition $\ker(H_{11})\subseteq \ker(H_{21})$ and by denoting $E=H_{21}+S_{21}$, 
   we obtain 
   \[
    \inf_{A\in \mathcal{N}_\succeq^n}\|AX-B\|_F^2\geq\inf_{\substack{H_{11}\succeq 0 \\ S_{11}=-S_{11}^\top}} \|(H_{11}+S_{11})\Sigma_1-U_1^\top BV_1\|_F^2+\inf_{E\in \R^{(n-r)\times r}}\|E\Sigma_1-U_2^\top BV_1\|_F^2+\|BV_2\|_F^2.
   \]
   Since the second infimum in the right hand side is $0$ and it is uniquely attained for $E_\star=H_{21}+S_{21}=U_2^\top BV_1\Sigma_1^{-1}:=Z$, it holds that 
   \[
    \inf_{A\in \mathcal{N}_\succeq^n}\|AX-B\|_F^2\geq\|\hat{A}_{11}\Sigma_1-U_1^\top BV_1\|_F^2+\|BV_2\|_F^2, 
   \]
   where $\hat{A}_{11}=\hat{H}_{11}+\hat{S}_{11}\in \mathcal{N}_\succeq^r$ is the unique solution of the first infimum (Lemma~\ref{lem:uniqueNSPSDP}), 
   and $\hat{H}_{11},\hat{S}_{11}$ are its symmetric and skew-symmetric part, respectively. Thus the conditions found are
   \begin{equation}
    \label{eq:condHS}
    H_{11}=\hat{H}_{11}, \quad  S_{11}=\hat{S}_{11}, \quad H_{21}+S_{21}=Z, \quad \ker(\hat{H}_{11})\subseteq\ker(H_{21}), \quad H_{22}-H_{21}\hat{H}_{11}^\dagger H_{21}^\top\succeq 0,
   \end{equation}
   together with the underlying constraint that $S_{22}$ is skew-symmetric. The restriction of $A\in U\mathcal{C}U^\top$ is introduced here and it implies that $H_{21}-S_{21}=0$. Thus all the constraints yields that the solution is of the form
   \[\reviseLAA{
    A_\star=U\left(
    \left(\begin{matrix}
     \hat{H}_{11} & \frac{1}{2}Z^\top \\
     \frac{1}{2}Z & K \\
    \end{matrix}\right)+
    \left(\begin{matrix}
     \hat{S}_{11} & -\frac{1}{2}Z^\top \\
     \frac{1}{2}Z & R \\
    \end{matrix}\right)\right)U^\top=U
    \left(\begin{matrix}
     \hat{A}_{11} & 0 \\
     Z & K+R \\
    \end{matrix}\right)U^\top,} 
   \]
   \noindent where $K \succeq \frac{1}{4}ZH_{11}^\dagger Z^\top$ and $R\in \R^{(n-r)\times (n-r)}$ is skew-symmetric. The imposition of the constraint $A\in U\mathcal{C}U^\top$ implies that the infimum of the problem is not always attained. Indeed the condition of \eqref{eq:condHS} that $\ker(\hat{H}_{11})\subseteq \ker(Z)$ is not fulfilled in general, especially if $\ker(\hat{H}_{11})\neq \{0\}$. For the continuation of the proof, we refer to \cite[Theorem 2]{baghel2022non}, where all the sub-cases (infimum attained or not) are discussed in detail.
  \end{proof}
 \end{thm}
 
 \begin{rem} 
  \label{rem:largenorm}
  Besides solving a different NSPSDP problem with the additional constraint $A\in \mathcal{N}_\succeq^{U\mathcal{C}U^\top}$, the issue of the method by Baghel et al.\ concerns the norm of the solution found when the infimum of \eqref{prob:preNSPSDP} is not attained, that is, when $X$ has low-rank and  $\ker(\hat{H}_{11})\neq \{0\}$ is not contained in $\ker(Z)$. 
  In that case, the matrix $K^{(\varepsilon)}\succeq \frac{1}{4}Z(\hat{H}^{(\varepsilon)}_{11})^{-1} Z^\top$ used to define $A_\star^{(\varepsilon)}$ has a very large norm, going to infinity as $\varepsilon$ goes to zero. 
  This leads to a huge norm solution, as it will be shown in the numerical examples. 
 \end{rem}

 \section{A new approach for the NSPSDP problem}
 \label{sec:newappr}
 
 In Section~\ref{sec:overview}, we have detected that an additional constraint $A\in U\mathcal{C}U^\top$ had been introduced in the optimization problem studied in \cite{baghel2022non}, changing it from \eqref{prob:NSPSDP} to \eqref{prob:preNSPSDP}. Thus in the following result, we fix this issue and uncover the remarkable property that,  as opposed to \eqref{prob:PSDP}, 
 \begin{center}
     \emph{the infimum of the NSPSDP problem~\eqref{prob:NSPSDP} is always attained.} 
 \end{center} 
 We characterize the set of minimizers of the problem and we determine the family of the solutions for the original NSPSDP problem.
 
 \begin{thm}
  \label{thm:solNSPSDP}
  Let $X,B\in \R^{n\times m}$ and consider the SVD 
  \[
   X=U\Sigma V^\top:=
   \begin{pmatrix}
    U_1 & U_2 \\
   \end{pmatrix}
   \begin{pmatrix}
    \Sigma_1 & 0 \\
    0 & 0 \\
   \end{pmatrix}
   \begin{pmatrix}
    V_1 & V_2 \\
   \end{pmatrix}^\top,
  \]
  where $\Sigma_1\in \R^{r\times r}$ has full rank. Then
  \[
   \inf_{A \in \mathcal{N}_\succeq^n} \|A X-B\|^2_F=\min_{A_{11}\in \mathcal{N}_\succeq^r}\|A_{11}\Sigma_1-U_1^\top B V_1\|_F^2+\|BV_2\|_F^2.
  \]
  Moreover, let us introduce
  \[
   \hat{A}_{11}=\argmin_{A_{11}\in \mathcal{N}_\succeq^r}\|A_{11}\Sigma_1-U_1^\top B V_1\|_F^2:=\hat{H}_{11}+\hat{S}_{11},
  \]
  where $\hat{H}_{11}$ and $\hat{S}_{11}$ are, respectively, the symmetric and skew-symmetric part of $\hat{A}_{11}$. Then any solution of the problem \eqref{prob:NSPSDP} is of the form
  \[
   A_\star=U
   \begin{pmatrix}
    \hat{A}_{11} & M^\top \\
    Z & N+G \\
   \end{pmatrix}U^\top=U_1\hat{A}_{11}U_1^\top+U_2ZU_1^\top+U_1M^\top U_2^\top+U_2(N+G)U_2^\top,
  \]
  where $Z=U_2^\top B V_1\Sigma_1^{-1}\in \R^{(n-r)\times r}$, 
  $M\in \R^{(n-r)\times r}$, $G\in \R^{(n-r)\times (n-r)}$ is skew-symmetric, $N\in \R^{(n-r)\times (n-r)}$ is symmetric, and they fulfil the conditions
  \[
   \ker(\hat{H}_{11})\subseteq \ker(Z+M), \qquad N \succeq \frac{1}{4}(Z+M)H_{11}^\dagger (Z+M)^\top,
  \]
  \reviseLAA{so that $A_\star\in \mathcal{N}_\succeq^n$.}
  \begin{proof}
   We proceed as in Theorem~\ref{thm:preNSPSDP}, until the introduction of the unnecessary constraint that the solution belongs to $U\mathcal{C}U^\top$. With the same notation, let $A\in \mathcal{N}_\succeq^n$ and set
   \[
    \hat{A}:=U^\top AU=U(H+S)U^\top:=U\left(
    \begin{pmatrix}
     H_{11} & H_{21}^\top \\
     H_{21} & H_{22} \\
    \end{pmatrix}+
    \begin{pmatrix}
     S_{11} & -S_{21}^\top \\
     S_{21} & S_{22} \\
    \end{pmatrix}\right) U^\top.
   \] 
   Let $\hat{A}_{11}\in \mathcal{N}_\succeq^r$ be the unique solution (provided by Lemma \ref{lem:uniqueNSPSDP}) of
   \[
    \inf_{A\in \mathcal{N}_\succeq^n}\|AX-B\|_F^2=\|\hat{A}_{11}\Sigma_1-U_1^\top BV_1\|_F^2+\|BV_2\|_F^2, 
   \]
   and denote by $\hat{H}_{11},\hat{S}_{11}$ its symmetric and skew-symmetric part, respectively. The conditions on the sub-matrices of the solution found in the proof of Theorem~\ref{thm:preNSPSDP} are
   \[
    H_{11}=\hat{H}_{11}, \qquad  S_{11}=\hat{S}_{11}, \qquad H_{21}+S_{21}=U_2^\top BV_1\Sigma_1^{-1}=Z, \qquad \ker(\hat{H}_{11})\subseteq \ker(H_{21}).
   \]
   In particular it is not required that $H_{21}=S_{21}=\frac{1}{2}Z$, but only that their sum is equal to $Z$.
   Thus this constraint simply implies that $A_\star$ is of the form
   \begin{equation}
    \label{eq:solformNSPSDP}
    A_\star=U
    \left(\begin{matrix}
     \hat{A}_{11} & M^\top \\
     Z & N+G \\
    \end{matrix}\right)U^\top=U_1\hat{A}_{11}U_1^\top+U_2ZU_1^\top+U_1M^\top U_2^\top+U_2(N+G)U_2^\top,
   \end{equation}
   where $M\in \R^{(n-r)\times r}$, $G\in \R^{(n-r)\times (n-r)}$ is skew-symmetric, $N\in \R^{(n-r)\times (n-r)}$ is symmetric and they fulfil the conditions
   \begin{equation}
    \label{eq:solconNSPSDP}
    \ker(\hat{H}_{11})\subseteq \ker(Z+M), \qquad N \succeq \frac{1}{4}(Z+M)H_{11}^\dagger (Z+M)^\top.
   \end{equation}
   This ensures that all $A_\star$ of the form of \eqref{eq:solformNSPSDP} belong to $\mathcal{N}_\succeq^n$. 
   Since they satisfy \eqref{eq:minpreNSPSDP}, they are all solutions of the optimization problem \eqref{prob:NSPSDP}. 
  \end{proof}
 \end{thm}
 
 Theorem~\ref{thm:solNSPSDP} shows that 
 \begin{equation} \label{eq:setFcaract}
  U\mathcal{F}U^\top:=\{UFU^\top: F\in \mathcal{F}\}=\argmin_{A \in \mathcal{N}_\succeq^n} \|A X-B\|^2_F,
 \end{equation}
 where $\mathcal{F}$ is the matrix family defined as
 \begin{align*}
       \mathcal{F} = \Bigg\{
  \left(\begin{matrix}
   \hat{A}_{11} & M^\top \\
   Z & N+G \\
  \end{matrix}\right) : 
  & \; G=-G^\top, 
  N \succeq \frac{1}{4}(Z+M)H_{11}^\dagger (Z+M)^\top, \\ 
  & \; 
    Z=U_2^\top B V_1\Sigma_1^{-1}, \ker(\hat{H}_{11})\subseteq \ker(Z+M) 
   \; \; \Bigg\}. 
 \end{align*} 
 We notice that $\mathcal{F}$ is not empty, since for instance
 \[
  F:=\begin{pmatrix}
     \hat{A}_{11} & -Z^\top \\
     Z & 0 \\
    \end{pmatrix}\in \mathcal{F}, 
 \]
by choosing $M=-Z$ and $N=0$ that fulfil the constraints given in \eqref{eq:solconNSPSDP}. 
This provides a solution $A=UFU^\top$ of \eqref{prob:NSPSDP} that attains the infimum, but $A\notin U\mathcal{C}U^\top$ and hence it cannot be found by the method proposed in \cite{baghel2022non}. 
 
 Now that we have introduced a new family of solutions for the NSPSDP problem \eqref{prob:NSPSDP}, we provide a simplified characterization of the matrix family $\mathcal{F}$, in which the constraint concerning the kernels in \eqref{eq:solconNSPSDP} is removed. In this way it will be also easier to highlight low-norm and low-rank properties of matrices in the family. 
 \begin{lem}
  \label{lem:familyNSPSDP}  
  With the same notation as in Theorem~\ref{thm:solNSPSDP}, the set $\mathcal{F}$ that can be used to characterize optimal solutions of~\eqref{prob:NSPSDP} \reviseLAA{(see~\eqref{eq:setFcaract})} is given by 
    \[
   \mathcal{F}=\left\{
   \begin{pmatrix}
    \hat{A}_{11} & (YW_1^\top-Z)^\top \\
    Z & N+G \\
   \end{pmatrix} : 
   G=-G^\top, 
   N \succeq \frac{1}{4}Y\Lambda^{-1} Y^\top, 
    Y\in \R^{(n-r)\times s}, 
    Z = U_2^\top B V_1\Sigma_1^{-1} 
   \right\}, 
  \]
  where 
  \[
   \frac{\hat{A}_{11}+\hat{A}_{11}^\top}{2}=\hat{H}_{11}:=
   \begin{pmatrix}
    W_1 & W_2 \\
   \end{pmatrix}
   \begin{pmatrix}
    \Lambda & 0\\
    0 & 0 \\
   \end{pmatrix}
   \begin{pmatrix}
    W_1 & W_2 \\
   \end{pmatrix}^\top
  \]
  is an SVD of $\hat{H}_{11}$, $s\leq r$, $\Lambda\in \R^{s\times s}$ is invertible, $W_1\in \R^{r\times s}$ and $W_2\in \R^{r\times (r-s)}$.  
  \begin{proof}
   The condition on the null spaces of $\hat{H}_{11}$ and $Z+M$ in \eqref{eq:solconNSPSDP} is equivalent to state that there exists a matrix $Y\in \R^{(n-r)\times s}$ such that $Z+M=YW_1^\top$. Indeed $W:=\begin{pmatrix} W_1 & W_2\\ \end{pmatrix}\in \R^{r\times r}$ is an orthonormal matrix and $\ker(\hat{H}_{11})=W_2$ implies
   \[
    (Z+M)W_2=YW_1^\top W_2=0.
   \]
   Thus
   \[
    N\succeq \frac{1}{4}YW_1^\top(W_1 \Lambda^{-1} W_1^\top)W_1 Y^\top=\frac{1}{4}Y\Lambda^{-1} Y^\top.
   \]   
  \end{proof}
 \end{lem}

 \subsection{Minimum-norm solution for the NSPSDP problem}

 Now we look for the matrices in the set of optimal solutions of the NSPSDP problem, namely $U\mathcal{F}U^\top$, with the smallest possible Frobenius norm. That is, we want to minimize the Frobenius norm of 
 \[
  A=U\begin{pmatrix}
    \hat{A}_{11} & (YW_1^\top-Z)^\top \\
    Z & N+G \\
   \end{pmatrix}U^\top\in U\mathcal{F} U^\top 
 \]
 \reviseLAA{where all the matrices fulfil the conditions in Lemma \ref{lem:familyNSPSDP}. Let us introduce the differentiable function
 \begin{equation}
  \label{eq:deffY}
  f(Y) := 
  \underbrace{\|Y W_1^\top-Z\|_F^2}_{=: g(Y)} + \frac{1}{16} \underbrace{\left\|Y \Lambda^{-1} Y^\top\right\|_F^2}_{=: h(Y)}, \qquad Y \in \R^{(n-r)\times s}.
 \end{equation}
 Note that $f$ implicitly depends on the fixed matrices $X$ and $B$, since $W_1$, $Z$ and $\Lambda$ do.
 
 \begin{lem}
  \label{lem:convf}
  The function $f$ defined in \eqref{eq:deffY} is strictly convex in $Y$ and its minimum is unique.
  \begin{proof}
   We first show that $g$ is strictly convex and that $h$ is convex, therefore their linear combination $f$ is strictly convex. Then we prove that $f$ attains a minimum, which is unique by strict convexity.  
   
   Let us define $g_1(x)=x^2$ and $g_2(Y)=\|YW_1^\top-Z\|_F$. The function $g_1$ is strictly convex and increasing on $[0,+\infty)$, while the function $g_2$ is convex on $\R^{(n-r)\times s}$, since for any $\alpha\in(0,1)$ and any $Y_1,Y_2\in \R^{(n-r)\times s}$, 
   \[
    g_2(\alpha Y_1+(1-\alpha)Y_2)=\|\alpha(Y_1W_1^\top-Z)+(1-\alpha)(Y_2W_1^\top-Z)\|_F\leq \alpha g_2(Y_1)+(1-\alpha)g_2(Y_2).
   \]
   Thus $g=g_1\circ g_2$ is strictly convex. In fact, its Hessian matrix is the identity, meaning that it is perfectly conditioned. The unique minimizer of $g(Y)$ is given by $ZW_1$. 
   
   Let us now prove the convexity of the quartic function $h(Y)=\|Y\Lambda^{-1}Y^\top\|_F^2=\trace((Y\Lambda^{-1}Y^\top)^2)$. Define the diagonal matrix $D:=\Lambda^{-1}$, which has, by construction, positive diagonal entries, and let
   \[
    h_1(X)=\trace(X^2), \qquad h_2(Y)=YDY^\top.
   \]
   The function $h_1$ is operator monotone (see~\cite[Definition 2.1]{carlen2010trace}), that is, for any $X_1,X_2 \in \mathcal{S}_\succeq^{n-r}$ such that $X_1\succ X_2$, we have $h_1(X_1)>h_1(X_2)$. Indeed, if $\lambda_i$ and $\mu_i$ are the eigenvalues of $X_1$ and $X_2$, respectively, 
   \[
    \trace(X_1^2-X_2^2)=\sum_{i=1}^{n-r} \lambda_i^2-\mu_i^2=\sum_{i=1}^{n-r} (\lambda_i+\mu_i)(\lambda_i-\mu_i)\geq (n-r)\min_i (\lambda_i+\mu_i)\trace(X_1-X_2)>0
   \]
   since $X_1-X_2\succ 0$. Moreover $h_1$ is a convex function, since it fulfills the hypothesis of \cite[Theorem 2.10]{carlen2010trace}, which states that if $\varphi:\R\rightarrow \R$ is convex, then  $X\rightarrow\trace(\varphi(X))$ is convex, where $\varphi(X)$ is meant as a matrix function; in this case $\varphi(x)=x^2$. The function $h_2$ is operator convex, that is, for any $\alpha\in(0,1)$ and any $Y_1,Y_2\in \R^{(n-r)\times s}$, 
   \[
    \alpha g_2(Y_1)+(1-\alpha) g_2(Y_2)-g_2(\alpha Y_1-(1-\alpha) Y_2) 
   \]
   \[
    =\alpha(1-\alpha)\left(Y_1D Y_1^\top+Y_2DY_2^\top-Y_1DY_2^\top-Y_2DY_1^\top\right)=\alpha(1-\alpha)(Y_1-Y_2)D(Y_1-Y_2)^\top\succeq 0.
   \] 
   Thus, for any $\alpha\in(0,1)$ and any $Y_1,Y_2\in \R^{(n-r)\times s}$, the properties of $h_1$ and $h_2$ yield the convexity of $h$:
   \[
    h(\alpha Y_1+(1-\alpha)Y_2)= h_1(h_2(\alpha Y_1+(1-\alpha)Y_2))\leq h_1(\alpha h_2(Y_1)+(1-\alpha)h_2(Y_2))\leq \alpha h(Y_1)+(1-\alpha)h(Y_2).
   \]
   
   The function $f$ is trivially bounded from below by $0$, and it is continuous with bounded level sets, since for any $\rho>0$
   \[
    f(Y)\leq \rho \quad \Rightarrow \quad \|Y\|_F=\|YW_1^\top\|_F\leq \|Y W_1^\top-Z\|_F+\|Z\|_F \leq \sqrt{\rho}+\|Z\|_F. 
   \] 
   Thus its infimum is attained while strict convexity implies its uniqueness. 
  \end{proof}
 \end{lem}

 Thanks to the properties of $f(Y)$ shown in Lemma \ref{lem:convf}, the following result provides the unique minimum Frobenius norm solution of  \eqref{prob:NSPSDP}, and it also shows its low-rank properties.}
 
 \begin{thm}
  \label{thm:minnormNSPSDP} 
  Using the same notation as in Theorem~\ref{thm:solNSPSDP}, Lemma~\ref{lem:familyNSPSDP} and Lemma~\ref{lem:convf}, the minimum Frobenius norm solution of \eqref{prob:NSPSDP} is 
  \[
   A_\star:=U
   \begin{pmatrix}
    \hat{A}_{11} & (Y_\star W_1^\top-Z)^\top \\
    Z & \frac{1}{4}Y_\star \Lambda^{-1} Y_\star^\top  \\
   \end{pmatrix}U^\top,
  \]
  where  $Y_\star$ is the unique global minimizer of 
  \[
   \min_{Y\in \R^{(n-r)\times s}} \|Y W_1^\top-Z\|_F^2+\frac{1}{16}\|Y \Lambda^{-1} Y^\top\|_F^2 = \min_{Y\in \R^{(n-r)\times s}} f(Y).
  \] 
  Moreover the symmetric part of $A_\star$ has rank $s$, which is the smallest possible rank of the symmetric part of a matrix in $U\mathcal{F}U^\top$, while the skew-symmetric part of $A_\star$ has rank not greater than $2r$.
  \begin{proof} 
   Since the Frobenius norm is unitary invariant and by using the characterization of Lemma~\ref{lem:familyNSPSDP}, we obtain 
   \[
    \|A_\star\|_F^2=\|\hat{A}_{11}\|_F^2+\|Z\|_F^2+\|Y_\star W_1^\top-Z\|_F^2+\|N+G\|_F^2.
   \]
   The first two addends of the sum are fixed. For the fourth, it holds that 
   \[
    \|N+G\|_F^2=\|N\|_F^2+2\textnormal{trace}(NG)+\|G\|_F^2=\|N\|_F^2+\|G\|_F^2,
   \]
   since $\textnormal{trace}(NG)=-\textnormal{trace}(NG)=0$ and thus the choice $G=0$ is optimal. Regarding the choice of $N$, recall that the definition of $\mathcal{F}$ implies that $N=\Delta+\frac{1}{4}Y\Lambda^{-1}Y^\top$ for some $\Delta\succeq 0$. As shown in \cite[Lemma 3]{baghel2022non}, the unique solution of the minimization problem
   \[
    \min_{\Delta\in \mathcal{N}_\succeq^n} \left\|\Delta+\frac{1}{4}Y \Lambda^{-1} Y^\top\right\|_F^2
   \]
   is $\Delta=0$ and hence the optimal choice $N = \frac{1}{4}Y\Lambda^{-1}Y^\top$ yields
   \begin{equation}
    \label{eq:infY}
    \inf_{A\in U\mathcal{F}U^\top}\|A\|_F^2=\|\hat{A}_{11}\|_F^2+\|Z\|_F^2+\inf_{Y\in \R^{(n-r)\times s}} \left(\|Y W_1^\top-Z\|_F^2+\frac{1}{16}\left\|Y \Lambda^{-1} Y^\top\right\|_F^2\right). 
   \end{equation}
   \reviseLAA{Hence, minimizing the norm of $A\in U\mathcal{F}U^\top$ is equivalent, up to constant summands, to minimize $f(Y)$. Lemma \ref{lem:convf} shows that the minimum of $f$ is uniquely attained and so $A_\star$ is the unique minimizer of $\min_{A\in U\mathcal{F}U^\top} \|A\|_F^2$. }
   
   The symmetric part of $A_\star$ is
   \[
    H_\star=\frac{A_\star+A_\star^\top}{2}=
    \begin{pmatrix}
     \hat{H}_{11} & \frac{1}{2}(Y_\star W_1^\top)^\top \\
     \frac{1}{2}Y_\star W_1^\top & \frac{1}{4}Y_\star\Lambda^{-1}Y_\star^\top \\
    \end{pmatrix}=
    \begin{pmatrix}
     W_1 \\
     \frac{1}{2}Y_\star \\ 
    \end{pmatrix}
    \begin{pmatrix}
     \Lambda & I_s \\
     I_s & \Lambda^{-1} \\
    \end{pmatrix}
    \begin{pmatrix}
     W_1^\top & \frac{1}{2}Y_\star^\top \\ 
    \end{pmatrix},
   \]
   where $I_s$ denotes the $s\times s$ identity matrix, and hence 
   \[
    \rank(H_\star)=\rank 
    \begin{pmatrix}
     \Lambda & I_s \\
     I_s & \Lambda^{-1} \\
    \end{pmatrix} = s.
   \]
   Moreover any matrix in $\mathcal{F}$ has symmetric part with rank not smaller than $s$, since using the Schur complement yields
   \[
    \rank 
    \begin{pmatrix}
     \hat{H}_{11} & \frac{1}{2}(Y_\star W_1^\top)^\top \\
     \frac{1}{2}Y_\star W_1^\top & N \\
    \end{pmatrix} \geq \rank 
    \begin{pmatrix}
     \hat{H}_{11} & 0 \\
     0 & N-\frac{1}{4}Y_\star \Lambda^{-1} Y_\star^\top \\
    \end{pmatrix} \geq s.
   \]
   To conclude the proof, we have for the skew-symmetric part of $A_\star$ that 
   \[
    \rank\left(\frac{A_\star-A_\star^\top}{2}\right)=
    \rank
    \begin{pmatrix}
     \hat{S}_{11} & \frac{1}{2}(Y_\star W_1^\top-2Z)^\top \\
     \frac{1}{2}(2Z-Y_\star W_1^\top) & 0 \\
    \end{pmatrix} \leq 2r.
   \]
  \end{proof}
 \end{thm}

 \subsection{Computing and approximating the minimum-norm solution} 

 \reviseLAA{Theorem~\ref{thm:minnormNSPSDP} describes the minimum-norm and minimum-rank solution of \eqref{prob:NSPSDP}. However a closed-form for the unique solution (Lemma~\ref{lem:convf}) 
 of the problem, 
 \begin{equation}
  \label{prob_minnorm}
  Y_\star=\argmin_{Y\in \R^{(n-r)\times s}} \|Y W_1^\top-Z\|_F^2+\frac{1}{16}\|Y \Lambda^{-1} Y^\top\|_F^2=\argmin_{Y\in \R^{(n-r)\times s}} f(Y), 
 \end{equation}
 is not available. We propose two different approaches to overcome this issue.

 \paragraph{Computing the minimum-norm solution}  
 
 To solve~\eqref{prob_minnorm}, we use a standard  non-linear optimization algorithm, namely gradient descent. The gradient of $f$ can be computed (see, e.g.,  \cite{petersen2008matrix}) explicitly: 
 \[
  \nabla f(Y)=\frac{1}{4} Y\Lambda^{-1}Y^\top Y \Lambda^{-1}+2(YW_1^\top-Z)W_1=\frac{1}{4} Y\Lambda^{-1}Y^\top Y \Lambda^{-1}+2Y-2ZW_1.
 \] 
 Note that $f(Y)$ is not globally Lipschitz because of the quartic term $\|Y \Lambda^{-1} Y^\top\|_F^2$. Although gradient descent is a relatively naive method, it worked very well for our purpose, because the problem is typically well conditioned (the Hessian is the identity plus a semidefinite matrix). In particular a simple backtracking line search will work well in our case. We update the current approximation $Y$ by choosing a suitable step size ${\gamma}$ that satisfies the Armijo condition:  
 \begin{equation}
  \label{eq:Armijo}
  f(Y-\gamma \nabla f(Y))\leq f(Y)-c\gamma\|\nabla f(Y)\|^2,
 \end{equation}
 where $0<c<1$. This allows us to get global convergence towards $Y_\star$, since $f$ is trivially bounded from below by $0$ and $\nabla f$ is locally Lipschitz continuous (see \cite[Theorem 3.2]{nocedal1999numerical}). We will provide the details in Algorithm~\ref{alg:Ystar}.
 
 \paragraph{Replacing the minimum-norm solution using Cardano's formula}  
 
 Instead of approximating $Y_\star$ using gradient descent applied on~\eqref{prob_minnorm}, we propose  the following approximation: 
 \[
  Y=\alpha ZW_1, \qquad \alpha\in [0,1].
 \]}
 This choice represents a convex combination of the minimizer of the first term $\|YW_1^\top-Z\|_F^2$ in \eqref{prob_minnorm}, namely $Y = ZW_1$, and of the minimizer of the second term, namely $Y = 0$. 
 The following lemma provides the optimal choice for the parameter $\alpha$, which minimizes $f(\alpha ZW_1)$. 
 \begin{lem}
  \label{lem:Cardano}
  Let 
  \[
   \varphi(\alpha):=f(\alpha ZW_1)=\frac{\alpha^4}{16}\|Z W_1\Lambda^{-1}W_1^\top Z^\top\|_F^2+\alpha^2\|ZW_1W_1^\top\|_F^2-2\alpha \langle ZW_1W_1^\top,Z\rangle+\|Z\|_F^2.
  \] 
  Then, the unique minimum of $\varphi$ is
  \begin{equation}
   \label{eq:cardano}
   \alpha_\diamond:=\argmin_{\alpha\in \R} \varphi(\alpha)=\sqrt[3]{\frac{p}{2}+\sqrt{\frac{p^2}{4}+\frac{p^3}{27}}}+\sqrt[3]{\frac{p}{2}-\sqrt{\frac{p^2}{4}+\frac{p^3}{27}}}\in (0,1), 
  \end{equation}
  where
  \[
   p=\frac{8\|ZW_1\|_F^2}{\|Z W_1\Lambda^{-1}W_1^\top Z^\top\|_F^2}>0.
  \]
  \begin{proof}
   The best value for $\alpha$ can be computed by finding a real zero of the derivative of $\varphi$, given by 
   \[
    \varphi'(\alpha)=\frac{\alpha^3}{4}\|Z W_1\Lambda^{-1}W_1^\top Z^\top\|_F^2+2\alpha\|ZW_1\|_F^2-2\|ZW_1\|_F^2.
   \]
   Since $\varphi'(0)<0$ and $\varphi'(1)>0$, there must be a zero in $(0,1)$. The equation $\varphi'(\alpha)=0$ is equivalent to
   \begin{equation}
    \label{eq:cubic}
    \alpha^3+p\alpha-p=0,
   \end{equation}
   which is a cubic equation with discriminant $\Delta=-(4p^2+27p^3)$. Since $\Delta$ is negative, Cardano's formula (see, e.g., \cite{cardano2007rules}) implies that the unique solution of \eqref{eq:cubic} is $\alpha_\diamond$ defined in \eqref{eq:cardano}.
  \end{proof} 
 \end{lem}
  
 \begin{rem}
  We have observed numerically that the solutions $Y_\star$ and $Y_\diamond=\alpha_\diamond ZW_1$ are typically close. Moreover it is possible to give a simple bound for the difference of the functional evaluated in these two points: 
  \[
   \|ZW_1W_1^\top-Z\|_F^2\leq f(Y_\star) \leq f(Y_\diamond) \leq f(\alpha ZW_1), \qquad \forall \alpha\in \R
  \]
  and the choices $\alpha=0$ and $\alpha=1$ imply
  \[
   |f(Y_\star)-f(Y_\diamond)|\leq \min\left\{\frac{1}{16}\|ZW_1\Lambda^{-1}W_1^\top Z^\top\|_F^2,\|Z\|_F^2-\|ZW_1W_1^\top-Z\|_F^2\right\}.
  \]
  We will compare these two solutions in Section~\ref{sec:numexp}.  
 \end{rem}

 \section{An algorithm for the NSPSDP problem}
 \label{sec:algor}
 
 Relying on Lemma~\ref{lem:familyNSPSDP} and Theorem~\ref{thm:minnormNSPSDP}, we now propose an algorithm for the NSPSDP problem \eqref{prob:NSPSDP}. The method follows the approach from \cite{baghel2022non}, but it performs a different post-processing of the solution of the reduced $r \times r$ dimensional problem, in order to avoid the unnecessary constraint introduced in~\eqref{prob:preNSPSDP}, and hence avoid unbounded solutions when the solution of \eqref{prob:preNSPSDP} is not attained. Algorithm~\ref{alg:NSPSDP} describes our proposed algorithm, which attains small Frobenius norm solutions with low-rank properties.  
 \begin{algorithm}[ht!]
  \caption{Algorithm for the NSPSDP problem \eqref{prob:NSPSDP}}
  \label{alg:NSPSDP}
  \begin{description}
   \item[Input:] Two matrices $X,B\in \R^{n\times m}$.
   \item[Output:] An optimum minimum-norm solution $A_\star$ of \eqref{prob:NSPSDP}; see Theorem~\ref{thm:minnormNSPSDP}. 
  \end{description}
  \begin{algorithmic}[1]
  
   \State \label{step1} Compute a rank-revealing SVD of 
   $X=
   \begin{pmatrix}
    U_1 & U_2 \\
   \end{pmatrix}
   \begin{pmatrix}
    \Sigma_1 & 0 \\
    0 & 0 \\
   \end{pmatrix}
   \begin{pmatrix}
    V_1 & V_2 \\
   \end{pmatrix}^\top$, with $\Sigma_1 \in \mathbb{R}^{r \times r}$ where $r=\rank(X)$. 
   
   \State \label{step2} Solve the reduced $r\times r$ NSPSDP problem 
   \[
    \hat{A}_{11}=\argmin_{A_{11}\in \mathcal{N}_\succeq^r}\|A_{11}\Sigma_1-U_1^\top B V_1\|_F^2:=\hat{H}_{11}+\hat{S}_{11},
   \]
   by means of a fast gradient method~\cite{nesterov1983method}, exactly as done in~\cite{baghel2022non}.  
   
   \State \label{step3} Compute $Z=U_2^\top BV_1 \Sigma^{-1}$ and a thin SVD of $\hat{H}_{11}=W_1 \Lambda W_1^\top$, with 
   $\Lambda \in \mathbb{R}^{s \times s}$ where
   $s = \rank(\hat{H}_{11})$.
   
   \State \label{step4} Solve the optimization problem 
   \[
    Y_\star=\argmin_{Y\in \R^{(n-r)\times s}} \|Y W_1^\top-Z\|_F^2+\frac{1}{16}\|Y \Lambda^{-1} Y^\top\|_F^2, 
   \]
   using gradient descent with backtracking line search; see Algorithm~\ref{alg:Ystar}.  
   
   \State \label{step5} Return the matrix
   \[
    A_\star=
    \begin{pmatrix}
     U_1 & U_2
    \end{pmatrix}
    \begin{pmatrix}
     \hat{A}_{11} & (Y_\star W_1^\top-Z)^\top \\
     Z & \frac{1}{4}Y_\star \Lambda^{-1} Y_\star^\top  \\
    \end{pmatrix}
    \begin{pmatrix}
     U_1 & U_2
    \end{pmatrix}^\top.
   \]
  \end{algorithmic}
 \end{algorithm}
 
The first two steps of Algorithm~\ref{alg:NSPSDP} involving the solution of the reduced problem are the same as in~\cite{baghel2022non}.  
The novelty of our algorithm is the post-processing procedure described in Steps~\ref{step3}, \ref{step4} and~\ref{step5}, which ensures that the algorithm always provides a minimum of problem \eqref{prob:NSPSDP}, while \cite{baghel2022non} was only providing a large-norm approximate solution \reviseLAA{(see~\eqref{sol:unattained})}, 
in the case the infimum of~\eqref{prob:preNSPSDP} is not attained. 

\paragraph{Computation of $Y_\star$}

Step \ref{step4} in Algorithm \ref{alg:NSPSDP} is solved by means of a gradient descent procedure, as shown in Algorithm \ref{alg:Ystar}.

\begin{algorithm}[ht!]
  \caption{\reviseLAA{Backtracking line search with Armijo condition for solving problem \eqref{prob_minnorm}}}
  \label{alg:Ystar}
  \reviseLAA{
  \begin{description}
   \item[Input:] The matrices $Z\in \R^{(n-r)\times r}$, $W_1\in \R^{r\times s}$, $\Lambda\in \R^{s\times s}$, a stepsize parameter $\theta>1$, a parameter $c>0$, maximum number of iterations $i_{\max}$, maximum number of checks for the Armijo condition $j_{\max}$ and an error tolerance $\tau$ 
   \item[Output:] A solution $Y_\star$ of \eqref{prob_minnorm}, that minimizes   $f(Y)=\|Y W_1^\top-Z\|_F^2+\frac{1}{16}\left\|Y \Lambda^{-1} Y^\top\right\|_F^2$. 
  \end{description}
  \begin{algorithmic}[1]
   \State Initialize $Y_0=\alpha_\diamond ZW_1$ defined by the Cardano approximation \eqref{eq:cardano}. 
   
   \State Set $i=1$ and $e_0=1$; compute $f(Y_0)$ and
   \[
    \nabla f(Y_0)=\frac{1}{4} Y_0\Lambda^{-1}Y_0^\top Y_0 \Lambda^{-1}+2Y_0-2ZW_1.
   \]
   \While{$i<i_{\max}$ \textbf{and} $\|\nabla f(Y_i)\|_F\leq\tau$ \textbf{and} $e_i\leq 10^{-15}$}
   \State Set 
   \[
    \gamma=\frac{1}{10}\cdot\frac{\|Y_i\|_F}{\|\nabla f(Y_i)\|_F}.
   \]

   \For{$j=1,\dots,j_{\max}$}
   \State Compute $\widehat{Y}=Y_i-\gamma \nabla f(Y_i)$.
   \If {$f(\widehat{Y})\leq f(Y_i)-c\gamma\|\nabla f(Y_i)\|_F^2$}
   \State \textbf{break}
   \Else
   \State $\gamma=\gamma/\theta$.
   \EndIf
   \EndFor
   \State Update $Y_{i+1}:=\hat{Y}$ and compute $\nabla f(Y_{i+1})$.
   \State Set $i=i+1$ and $e_i=\frac{|f(Y_{i+1})-f(Y_i)|}{f(Y_i)}$.
   \EndWhile
   \State \textbf{return} $Y_\star:=Y_{i}$.
  \end{algorithmic}}
 \end{algorithm}
 \reviseLAA{The parameter $\theta$ introduced in the algorithm determines the reduction rate for the stepsize $\gamma$, the constant $c$ is used in the Armijo condition \eqref{eq:Armijo}, and the safety parameter $j_{\max}=50$ stops the algorithm in the  case the line search fails to fulfill \eqref{eq:Armijo}. In our experiments, we have always used $\theta=1.5$ and $c=10^{-4}$}. 
 In Section~\ref{sec:numexp}, we give more details about the stopping criterion for the algorithm.

\reviseLAA{\paragraph{Computing $Y_\diamond$}}  In Step \ref{step4} of Algorithm \ref{alg:NSPSDP} we may replace $Y_\star$ with $Y_\diamond = \alpha_\diamond ZW_1$, where $\alpha_\diamond$ is defined in \eqref{eq:cardano}. This approximate solution allows us to avoid solving~\eqref{prob_minnorm}, 
while still preserving the low-rank properties of $A_\star$, since 
 \[
  A_\diamond:=U
  \begin{pmatrix}
   \hat{A}_{11} & (Y_\diamond W_1^\top-Z)^\top \\
   Z & \frac{1}{4}Y_\diamond \Lambda^{-1} Y_\diamond^\top  \\
  \end{pmatrix}U^\top
 \]
 has a symmetric part of rank $s$ and skew-symmetric part of rank bounded by $2r$, as shown in the proof of Theorem~\ref{thm:minnormNSPSDP} for $Y_\star$. However the solution $A_\diamond$ is not guaranteed to have minimum Frobenius norm. In the numerical experiments in Section~\ref{sec:numexp}, we will compare these two approaches, both in terms of computational time and size of the norm of the solution.
 
 \subsection{Computational complexity}
 
 Given $X,B\in \R^{n\times m}$ with $\rank(X)=r$, let us analyze the computational complexity of Algorithm \ref{alg:NSPSDP}: 
 \begin{enumerate}
  \item The SVD in Step \ref{step1} requires $\io(mn\min(m,n))$ flops. 
  
  \item As shown in \cite{baghel2022non}, the complexity of the fast gradient method applied to the $r\times r$ NSPSDP problem is dominated by the computation of the gradient $AXX^\top-BX^\top$, which is $\io(\tau_1 r^3)$, where $\tau_1$ is the number of iterations performed.  See Section~\ref{sec:stopcrit} below for a discussion on the stopping criterion used. 
  
  \item For Step \ref{step3}, the cost of the computation of $Z$ is $\io(nmr)$, while the computation of the SVD of $\hat{H}_{11}$ is $\io(r^3)$.
 
  \item The computation of the solution of the optimization problem in Step \ref{step4} is performed by means of a gradient descent method, which requires the computation of the matrix $\nabla f(Y)=Y\Lambda^{-1}(Y^\top Y)\Lambda^{-1}+2(Y-ZW_1)$ at each iteration, leading to a cost of $\io(\tau_2 ((n-r)^2 r+s^3))$, where $\tau_2$ is the number of iterations performed. See Section~\ref{sec:stopcrit} below for a discussion on the stopping criterion used.  
  
  \item Step \ref{step5} requires matrix multiplications whose cost sums up to $\io(n^2r)$ flops.
  
 \end{enumerate}
 
 Thus the overall cost of the algorithm is $\io(mn\min(m,n) + \tau_2(n-r)^2 r + \tau_1 r^3)$. This cost reduces to $\io(mn\min(m,n) + \tau_1 r^3)$ if the solution of the optimization problem in Step 4 of 
 Algorithm~\ref{alg:NSPSDP} is replaced by $Y_\diamond=\alpha_\diamond ZW_1$, with 
 $\alpha_\diamond$ defined as in \eqref{eq:cardano}. 

Typically, $m \gg n \geq r$ since the number of variables is smaller than the number of measurements, and hence the most expensive step is the SVD computation in Step~\ref{step1} in $\io(mn\min(m,n))$ operations. 
We will provide numerical experiments on the scalability of Algorithm~\ref{alg:NSPSDP} in Section~\ref{sec:scalability}.

 \subsection{The algorithm for nearly low-rank $X$}
 
 The solution of the Procrustes problem \eqref{prob:NSPSDP} may be also modified in order to take advantage of a possible low-rank property of the input data. When the matrix $X$ is close to be rank-deficient, it is possible to replace it by its $(k-1)$-rank approximation $\widetilde{X}$ such that 
 \[
  \sigma_k\leq \eta \sigma_1,
 \]
 where $\sigma_i$ are the singular values of $X$ and $0<\eta\ll 1$ is a small tolerance. The following result gives an idea of how this modification changes the error $\|AX-B\|_F$ of the problem.
 
 \begin{thm}
  \label{thm:quasilowrank}
  Given $\eta>0$ , let $X=\widetilde{X}+\Delta X$ be a decomposition of the form
  \[
   X=
   U\begin{pmatrix}
    \sigma_1 & & &\\
    & \ddots & & \\
    & & \sigma_{k-1} & \\
    & & & 0 & & \\
    & & & & \ddots & & \\
    & & & & & 0
   \end{pmatrix}V^\top+
   U\begin{pmatrix}
    0 & & &\\
    & \ddots & & \\
    & & 0 & \\
    & & & \sigma_k & & \\
    & & & & \ddots & & \\
    & & & & & \sigma_q
   \end{pmatrix}V^\top , 
  \]
  where $q=\min(m,n)$ and $\sigma_k\leq \eta \sigma_1$. Define
  \reviseLAA{\[
   \tilde{\gamma}:=\min_{A\in \mathcal{N}^n_\succeq}\|A\tilde{X}-B\|_F:=\|\tilde{A}\tilde{X}-B\|_F, \qquad
   \gamma_\star:=\min_{A\in \mathcal{N}^n_\succeq}\|AX-B\|_F:=\|A_\star X-B\|_F,
  \]}
  where $A_\star$ is the unique minimizer of the original NSPSDP problem, while $\tilde{A}$ is a minimizer of the perturbed NSPSDP problem with smallest norm. Then
  \begin{equation}
   \label{eq:boundlr}
   \tilde{\gamma}-\xi \|A_\star\|_F\leq \gamma_\star \leq \tilde{\gamma}+\xi \|\tilde{A}\|_F,
  \end{equation}
  where $\xi=\eta \sigma_1 \sqrt{q-k+1}$.
  \begin{proof}
   The norm of the perturbation $\Delta X$ is bounded by $\xi$ since
   \[
    \|\Delta X\|_F=\sqrt{\sum_{i=k}^q \sigma_i^2} \leq \sqrt{\eta^2 (q-k+1) \sigma_1^2} \leq \eta \sigma_1 \sqrt{q-k+1}=\xi.
   \]
   Thus
   \[
    \gamma_\star=\|A_\star X-B\|_F \leq \|\tilde{A} (\tilde{X}+\Delta X)-B\|_F \leq \tilde{\gamma}+\|\tilde{A}\|_F\|\Delta X\|_F,
   \]
   \[
    \tilde{\gamma}=\|\tilde{A}\tilde{X}-B\|_F \leq \|A_\star (X-\Delta X)-B\|_F \leq \gamma_\star+\|A_\star\|_F\|\Delta X\|_F
   \]
   and the claim follows.
  \end{proof} 
 \end{thm}
 While the upper bound of \eqref{eq:boundlr} contains $\|\tilde{A}\|_F$, which is generally small, the lower bound provided depends on $\|A_\star\|_F$, which may be very large when $X$ has small non-zero singular values. Thus even a small perturbation can cause a huge change in the value of the functional $\|AX-B\|_F$. We will discuss more in detail this behaviour in the numerical experiments in Section~\ref{sec:numexp}.

 \section{Numerical experiments}
 \label{sec:numexp}

 In this section, we compare the performances of four different methods for solving the NSPSDP problem: 
 \begin{itemize}

  \item Fast Gradient Method (FGM), applied directly to~\eqref{prob:NSPSDP}, an optimal first-order method for convex problems developed in \cite{nesterov1983method}, as done in~\cite{baghel2022non}. 
 
  \item ANalytic Fast Gradient Method (ANFGM) uses the semi-analytical approach proposed in \cite{baghel2022non} to reduce the dimension of the problem from $n\times m$ to $r\times r$, solves it via FGM, and then uses the post-processing as in Theorem~\ref{thm:preNSPSDP}. This algorithm solves a more constrained NSPSDP problem, namely  \eqref{prob:preNSPSDP}.  

  \item MINimum Gradient Descent (MINGD) implements Algorithm~\ref{alg:NSPSDP}.  
  
  \item CARDano (CARD) follows the same approach as MINGD by implementing Algorithm~\ref{alg:NSPSDP}, but it replaced the minimum-norm solution in step~\ref{step4} by the solution obtained via the Cardano's formula; see \eqref{eq:cardano}. 
 \end{itemize}

 In~\cite{baghel2022non}, the performances of the interior point method (IPM) proposed in \cite{krislock2004local} and the convex \Matlab software (CVX) \cite{grant2008graph, grant2014cvx} was compared to FGM and ANFGM. 
 In most cases, IPM and CVX were slower and less accurate and hence we do not report their results here. The code and experiments are available from \url{https://github.com/StefanoSicilia/NS_Procr_min_norm}.

 Before providing a comparison on several types of data sets, let us discuss the stopping criterion and metrics used.

 \paragraph{Stopping criteria} \label{sec:stopcrit}
 
ANFGM, MINGD and CARD reduce the size of the original problem, and then rely on the fast gradient method (FGM) proposed in \cite{nesterov1983method}. 
 We will use the same stopping criterion as in \cite{baghel2022non}, namely when the $k$th iterate $A^{(k)}$ satisfies
 \[
  \|A^{(k)}-A^{(k-1)}\|_F<\delta \|A^{(1)}-A^{(0)}\|_F,
 \]
 with $\delta=10^{-6}$, that is, the modification of $A$ compared to the first step is less than $10^{-6}$, or when the maximum number of iteration $k_{\max}=10000$ is reached. We will use this stopping criterion for FGM, ANFGM, MINGD and CARD  in all numerical experiments.
 
 \reviseLAA{Regarding Algorithm \ref{alg:Ystar} for the computation of $Y_\star$ in MINGD, we stop the gradient descent algorithm when the iterate $Y^{(k)}$ satisfies one of the following conditions 
 \[
  \|\nabla(Y^{(k)})\|_F\leq \tau:= 10^{-8}, \qquad \frac{|f(Y^{(k+1)})-f(Y^{(k)})|}{f(Y^{(k)})}\leq 10^{-15}, 
 \]
 or when the number of iteration number exceeds $i_{\max}=500$. 
 The stopping criterion that is most often reached first is the one based on the norm of the gradient. However, in some cases, we observed that the objective function reached machine accuracy (about $10^{-15}$) while the norm of the gradient was significantly larger than $10^{-8}$. The reason is the presence of a quartic term in~$f$.}

\paragraph{Metrics used to compare the algorithms}
 
 We will use the the following metrics to compare the algorithms: given a solution $A$ computed by an algorithm, we will report 
 \begin{enumerate}
  \item The relative error 
  \[
   \frac{\|AX-B\|_F}{\|B\|_F}, 
  \]
   which measures the quality of a solution. 
  
  \item The Frobenius norm of the solution, $\|A\|_F$.
  
  \item The time required by the algorithm (in seconds).
  
  \item The rank of the symmetric part of the solution,  $\textnormal{rk}(\symm(A))$.
  
  \item The rank of the skew-symmetric part of the solution,  $\textnormal{rk}(\skeww(A))$.
 \end{enumerate}

 \subsection{Synthetic data}
 
 Similarly to \cite{baghel2022non, gillis2018semi}, we generate some synthetic data for different values of $n,m$. For each dimension, we generate 20 random matrices and report the mean and standard deviation of the results. 
Given $m$, $n$ and $r$, a rank $r$ matrix $X\in \R^{n\times m}$ is generated as the product of two randomly generated matrices of dimension $n \times r$ and $r \times n$ whose entries follow the standard normal distribution, that is,  
\texttt{X = randn(n,r)*randn(r,m)} in \Matlab, while the entries of $B$ are generated uniformly at random, \texttt{B = randn(n,m)} in \Matlab. 

We choose $r\ll \min(m,n)$ in order to highlight the low-rank properties of the algorithms and solutions. 
Tables \ref{tab:507020}, \ref{tab:10010040} and \ref{tab:20020050} show the results of the four algorithm to solve~\eqref{prob:NSPSDP} for various dimensions.  

\renewcommand{\arraystretch}{1.2}

 \begin{table}[ht!] 
 \centering
 \begin{tabular}{|c|c|c|c|c|c|}
 \hline
       & Relative error     & $\|A\|_F$                & Time (sec.)    & $\textnormal{rk}(\symm(A))$ & $\textnormal{rk}(\skeww(A))$ \\ \hline
 ANFGM & $0.8674 \pm0.0114$ & $(2.07 \pm 0.48)\cdot {\color{red}10^7}$ & $\mathbf{<0.01}$ & $  20 \pm   0$       & $40 \pm0$             \\ \hline
 FGM   & $0.8674 \pm0.0114$ & $1.2845 \pm      0.0646$ & $0.86 \pm0.44$ & $39.40 \pm0.68$       & $40 \pm0$             \\ \hline
 MINGD & $0.8674 \pm0.0114$ & \reviseLAA{$\mathbf{0.9116 \pm      0.0576}$} & $\mathbf{<0.01}$ & $\mathbf{9.40 \pm0.68}$       & $40 \pm0$            \\ \hline
 CARD  & $0.8674 \pm0.0114$ & $0.9408 \pm      0.0582$ & $\mathbf{<0.01}$ & $\mathbf{ 9.40 \pm0.68}$       & $40 \pm0$             \\ \hline
 \end{tabular}
 \caption{Results on a sample of 20 matrices, with $n=50$, $m=70$, and $r=20$. Best solution highlighted in bold, if applicable.}
 \label{tab:507020}
 \end{table}
 
 \begin{table}[ht!] 
 \centering
 \begin{tabular}{|c|c|c|c|c|c|}
 \hline
       & Relative error     & $\|A\|_F$                & Time (sec.)    & $\textnormal{rk}(\symm(A))$ & $\textnormal{rk}(\skeww(A))$ \\ \hline
 ANFGM & $0.8007 \pm0.0060$       & $(2.12 \pm0.22)\cdot {\color{red}10^7}$ & $\mathbf{0.03 \pm0.01}$       & $40 \pm 0$               & $80 \pm0$                 \\ \hline
FGM   & $0.8007 \pm0.0060$       & $1.5392 \pm0.0461$              & $ 5.37 \pm1.86$       & $78.90 \pm 0.64$          & $80 \pm0$                 \\ \hline
MINGD & $0.8007 \pm0.0060$       & \reviseLAA{$\mathbf{1.1459 \pm0.0285}$}              & \reviseLAA{$0.05 \pm 0.03$}       & $\mathbf{18.90 \pm 0.64}$          & $80 \pm0$                 \\ \hline
CARD  & $0.8007 \pm0.0060$       & $1.1960 \pm0.0408$               & $\mathbf{0.03 \pm0.00}$       & $\mathbf{18.90 \pm 0.64}$          & $80 \pm0$                 \\ \hline
 \end{tabular}
 \caption{Results on a sample of 20 matrices, with $n=100$, $m=100$, and $r=40$. Best solution highlighted in bold, if applicable.}
 \label{tab:10010040}
 \end{table}

 \begin{table}[ht!]
 \centering
 \begin{tabular}{|c|c|c|c|c|c|}
 \hline
         & Relative error    & $\|A\|_F$                 & Time (sec.)   & $\textnormal{rk}(\symm(A))$ & $\textnormal{rk}(\skeww(A))$ \\ \hline
 ANFGM & $0.8776 \pm0.0028$      & $(9.22 \pm4.84)\cdot {\color{red} 10^6}$ & $\mathbf{0.03 \pm 0.00}$          & $   50 \pm   0$          & $100 \pm0$                \\ \hline
FGM   & $0.8776 \pm0.0028$      & $ 1.235 \pm 0.0198$ & $16.67 \pm6.34$          & $171.60 \pm4.07$          & $100 \pm0$                \\ \hline
MINGD & $0.8776 \pm0.0028$      & \reviseLAA{$\mathbf{0.7859 \pm 0.0154}$} &  \reviseLAA{$0.12 \pm 0.08$}          & $\mathbf{23.95 \pm0.83}$          & $100 \pm0$                \\ \hline
CARD  & $0.8776 \pm0.0028$      & $ 0.8203 \pm 0.021$ & $\mathbf{0.03 \pm 0.01}$          & $\mathbf{23.95 \pm0.83}$          & $100 \pm0$                \\ \hline
 \end{tabular}
 \caption{Results on a sample of 20 matrices, with $n=200$, $m=200$, and $r=50$. Best solution highlighted in bold, if applicable.} 
 \label{tab:20020050}
 \end{table}
 
 For all the choices of $n,m$ and $r$, we observe a similar behaviour of the algorithms: 
 \begin{itemize}
 
  \item All the methods provide the same relative error, up to a difference of order $10^{-9}$. Although ANFGM does not solve the same problem, it constructs a solution using the parameter $\varepsilon = 10^{-8}$ in~\eqref{sol:unattained}, and hence the solution generated has a relative error very close to that of the other algorithms, \reviseLAA{even though MINGD and CARD are preferable since they attain a minimum.}  
  
  \item The norm of the solution provided by ANFGM is very large (of order $10^7$), while, in contrast, FGM, MINGD and CARD have solutions with norms of comparable size, although MINGD provides always the smallest value, as shown in  Theorem~\ref{thm:minnormNSPSDP}. 
  
  \item The CPU time required by FGM is significantly larger (approximately 50 times) than that required by the other algorithms, because it does not reduce the size of the problem. 
  
  \item The symmetric part of the solution provided by MINGD and CARD have the lowest rank, followed by ANFGM. 
  On the other hand, the symmetric part of the solution of FGM is almost full rank. 
  
  \item All the skew-symmetric parts of the solutions have the same rank equal to $2r$.
  
 \end{itemize}
 
 In summary the best methods \revise{on our set of synthetic data} are MINGD and CARD: they provide a solution with small norm (the one from CARD is generally \reviseLAA{$4\%$} larger than that of MINGD), the computation is fast (MINGD is a bit slower than CARD; see the next sections for more experiments) and they satisfy the low-rank properties of Theorem~\ref{thm:minnormNSPSDP}, which makes the storage of the solution more efficient. Instead FGM is quite slow, while ANFGM provides solutions with very large norm and neither of these two methods provides a solution whose symmetric part has the  smallest possible rank.

 \subsection{Synthetic data with perturbed low-rank $X$}
  
 In order to study the stability of the approaches, we add a full-rank perturbation of order $\eta=10^{-8}$ to the previous synthetic data and we observe the behaviour of the algorithms. The perturbed matrix \texttt{X = randn(n,r)*randn(r,m)+$\eta$*randn(n,m)} is now full rank and this guarantees a unique solution to \eqref{prob:NSPSDP}. We also consider two other versions of MINGD and CARD, that will be denoted by MINGD$_{\textnormal{lr}}$ and CARD$_{\textnormal{lr}}$ respectively, where the matrix $X$ is replaced by its rank-$k$ approximation and $k$ is chosen such that 
 \[
  \sigma_{k+1}(X) \leq \eta \sigma_1(X),
 \]
 and we test the results given by Theorem \ref{thm:quasilowrank}. 
 
 In this way, we analyze if it is possible to exploit the nearly low-rank properties of the problem, see the discussion around Theorem~\ref{thm:quasilowrank}. Tables \ref{tab:per507020}, \ref{tab:per10010040} and \ref{tab:per20020050} show the results of the six algorithms.  
 
 \begin{table}[ht!] 
 \centering
 \begin{tabular}{|c|c|c|c|c|c|}
 \hline
       & Relative error     & $\|A\|_F$                  & Time (sec.)      & $\textnormal{rk}(\symm(A))$ & $\textnormal{rk}(\skeww(A))$ \\ \hline
 ANFGM & $\mathbf{0.7544 \pm0.0135}$ & $(5.08 \pm0.21)\cdot {\color{red}10^8}$ & $0.20 \pm0.01$ & $ 43.10 \pm0.64$      & $50 \pm0$             \\ \hline
 FGM   & $0.8703 \pm0.0146$ & $ 1.2692 \pm       0.0486$ & $2.06 \pm0.04$ & $39.45 \pm0.69$      & $50 \pm0$             \\ \hline
 MINGD & $\mathbf{0.7544 \pm0.0135}$ & $(5.08 \pm0.21)\cdot {\color{red}10^8}$ & $0.22 \pm0.04$ & $ 43.10 \pm0.64$      & $50 \pm0$             \\ \hline
 CARD  & $\mathbf{0.7544 \pm0.0135}$ & $(5.08 \pm0.21)\cdot {\color{red}10^8}$ & $0.22 \pm0.03$ & $ 43.10 \pm0.64$      & $50 \pm0$             \\ \hline
 MINGD$_{\textnormal{lr}}$ & $0.8703 \pm0.0146$ & \reviseLAA{$\mathbf{0.8986 \pm       0.0410}$} & $\mathbf{0.01 \pm0.00}$ & $\mathbf{9.45 \pm0.69}$      & $\mathbf{40 \pm0}$             \\ \hline
 CARD$_{\textnormal{lr}}$  & $0.8703 \pm0.0146$ & $0.9272 \pm       0.0441$ & $\mathbf{0.01 \pm0.00}$ & $\mathbf{9.45 \pm0.69}$      & $\mathbf{40 \pm0}$             \\ \hline
 \end{tabular}
 \caption{Perturbed $X$: results on a sample of 20 matrices, with $n=50$, $m=70$, and $r=20$. Best solution highlighted in bold, if applicable.}
 \label{tab:per507020}
 \end{table}
 
 \begin{table}[ht!] 
 \centering
 \begin{tabular}{|c|c|c|c|c|c|}
 \hline
       & Relative error     & $\|A\|_F$                   & Time (sec.)    & $\textnormal{rk}(\symm(A))$ & $\textnormal{rk}(\skeww(A))$ \\ \hline
 ANFGM                     & $\mathbf{0.6506 \pm0.0134}$ & $(2.16\pm0.98)\cdot {\color{red} 10^9}$ & $ 0.70 \pm0.01$ & $93.60 \pm1.82$ & $100 \pm0$ \\ \hline
FGM                       & $0.8064 \pm0.0073$ & $1.5450 \pm0.0469$                  & $8.02 \pm0.51$ & $78.50 \pm0.83$ & $100 \pm0$ \\ \hline
MINGD                     & $\mathbf{0.6506 \pm0.0134}$ & $(2.16\pm0.98)\cdot {\color{red} 10^9}$ & $0.89 \pm0.06$ & $93.60 \pm1.82$ & $100 \pm0$ \\ \hline
CARD                      & $\mathbf{0.6506 \pm0.0134}$ & $(2.16\pm0.98)\cdot {\color{red} 10^9}$ & $ 0.80 \pm0.01$ & $93.60 \pm1.82$ & $100 \pm0$ \\ \hline
MINGD$_{\textnormal{lr}}$ & $0.8064 \pm0.0073$ & \reviseLAA{$\mathbf{1.1504 \pm 0.0350}$}                & \reviseLAA{$0.04 \pm0.02$}    & $\mathbf{18.50 \pm0.83}$ & $ \mathbf{80 \pm0}$ \\ \hline
CARD$_{\textnormal{lr}}$  & $0.8064 \pm0.0073$ & $1.1940 \pm 0.0449$                 & $\mathbf{0.03 \pm0.00}$    & $\mathbf{18.50 \pm0.83}$ & $ \mathbf{80 \pm0}$ \\ \hline
 \end{tabular}
 \caption{Perturbed $X$: results on a sample of 20 matrices, with $n=100$, $m=100$, and $r=40$. Best solution highlighted in bold, if applicable.}
 \label{tab:per10010040}
 \end{table}
 
 \begin{table}[ht!] 
 \centering
 \begin{tabular}{|c|c|c|c|c|c|}
 \hline
       & Relative error     & $\|A\|_F$                     & Time (sec.)     & $\textnormal{rk}(\symm(A))$ & $\textnormal{rk}(\skeww(A))$ \\ \hline
 ANFGM                     & $\mathbf{0.6401 \pm0.0098}$      & $(3.62 \pm 1.11)\cdot {\color{red}10^9}$ & $ 2.21 \pm0.05$          & $192.25 \pm1.52$         & $200 \pm0$                \\ \hline
FGM                       & $0.8767 \pm0.0031$      & $          1.2380 \pm         0.0196$ & $32.07 \pm0.33$          & $ 171.70 \pm2.83$         & $200 \pm0$                \\ \hline
MINGD                     & $\mathbf{0.6401 \pm0.0098}$      & $(3.62 \pm 1.11)\cdot {\color{red}10^9}$ & $  2.20 \pm0.06$          & $ 192.70 \pm1.13$         & $200 \pm0$                \\ \hline
CARD                      & $\mathbf{0.6401 \pm0.0098}$      & $(3.62 \pm 1.11)\cdot {\color{red}10^9}$ & $ 2.19 \pm0.06$          & $ 192.70 \pm0.98$         & $200 \pm0$                \\ \hline
MINGD$_{\textnormal{lr}}$ & $0.8767 \pm0.0031$      & \reviseLAA{$ \mathbf{0.7863 \pm 0.0169}$} & \reviseLAA{$ 0.10 \pm 0.06$}          & $ \mathbf{24.30 \pm 0.80}$         & $\mathbf{100 \pm0}$                \\ \hline
CARD$_{\textnormal{lr}}$  & $0.8767 \pm0.0031$      & $         0.8223 \pm         0.0215$ & $\mathbf{0.03 \pm 0.00}$          & $ \mathbf{24.30 \pm 0.80}$         & $\mathbf{100 \pm 0}$                \\ \hline
 \end{tabular}
 \caption{Perturbed $X$: results on a sample of 20 matrices, with $n=200$, $m=200$, and $r=50$ Best solution highlighted in bold, if applicable.}
 \label{tab:per20020050}
 \end{table}
 
 For all the choices of $n,m$ and $r$ we observe a similar behaviour of the algorithms: 
 \begin{itemize}
  \item ANFGM, MINGD and CARD provide the same result, by computing the unique solution of the perturbed NSPSDP problem. Thus relative error, norm of the solution and the rank of its symmetric and skew-symmetric part coincide. 
  
  \item MINGD$_{\textnormal{lr}}$ and CARD$_{\textnormal{lr}}$ compute the solution of the low-rank problem associated to the low-rank approximation of $X$, and also FGM seems to do so. The reason is that FGM is a first-order method, with linear convergence of $\io((1-\kappa)^t)$, where $t$ is the iteration count and $\kappa=\frac{\sigma_r(X)}{\sigma_1(X)}\in (0,1)$ is the ratio between the $r$th and first singular value of $X$.
  Since the last $n-r$ are very small, it cannot converge within the allotted number of iterations.  
  These three methods have the same relative error, which is approximately $30-40\%$ larger than the unperturbed solutions.  
  However FGM is much slower than MINGD$_{\textnormal{lr}}$ and CARD$_{\textnormal{lr}}$.
  
  \item The norm of the solutions of ANFGM, MINGD and CARD are very large, in contrast to the small norms of the other methods. The reason is that $X$ has very small non-zero singular values. To be able to scale up these small contribution to approximate $B$, the norm of $A$ in these directions must be of the order of the inverse of these singular values. 
  
  \item ANFGM, MINGD and CARD are slower than  MINGD$_{\textnormal{lr}}$ and CARD$_{\textnormal{lr}}$. This is because they cannot reduce the problem size much since $X$ has full rank, while MINGD$_{\textnormal{lr}}$ and CARD$_{\textnormal{lr}}$ work with a $r \times r$ subproblem.
  
  \item Only the solution computed by MINGD$_{\textnormal{lr}}$ and CARD$_{\textnormal{lr}}$ have low-rank properties.
  
 \end{itemize}
 
 These experiments show that the low-rank approximation improves the speed of the algorithm and provides low-norm solution, with low-rank properties, but at the price of increasing the relative error. This is unavoidable, since the unique solution of the full-rank problem has  large norm. 
 

 \subsection{A large low-rank example}
 
 Now we consider an example of large dimension, with $n=500$, $m=10000$ and $r=10$. We define
 \begin{equation}
  \label{eq:JH}
  J_n=\left(
  \begin{matrix}
   1 & \cdots & 1 \\
   0 & \ddots & \vdots \\
   \vdots & \ddots & 1 \\
   0 & \cdots & 0 \\
  \end{matrix}\right)\in \R^{n\times r}, 
  \qquad H=\left(
  \begin{matrix}
   1 & 0 & \cdots & 0 \\
   2 & 1 & \ddots & \vdots \\
   \vdots & \ddots & \ddots & 0 \\
   n & \cdots & 2 & 1 \\
  \end{matrix}\right)\in \R^{n\times m}, 
 \end{equation}
 and we solve the NSPSDP problem \eqref{prob:NSPSDP} for $X=J_n J_m^\top$ and $B=H$. The results are shown in Table \ref{tab:JH}. 
 
 \begin{table}[ht!] 
 \centering
 \begin{tabular}{|c|c|c|c|c|c|}
 \hline
       & Relative error & $\|A\|_F$  & Time (sec.) & $\textnormal{rk}(\symm(A))$ & $\textnormal{rk}(\skeww(A))$ \\ \hline
 ANFGM & 0.9605         & $5.8633\cdot {\color{red}10^{15}}$ & 4.0469      & 500                         & 498                          \\ \hline
 FGM   & 0.9605         & $1.0172\cdot 10^4$ & 207.9887    & 321                         & 20                           \\ \hline
 MINGD & 0.9605         & $\mathbf{8.8618\cdot 10^3}$ & 3.8801      & $\mathbf{5}$                           & $\mathbf{12}$                          \\ \hline
 CARD  & 0.9605         & $8.8633\cdot 10^3$ & \textbf{3.6061}      & $\mathbf{5}$                           & $\mathbf{12}$                           \\ \hline
 \end{tabular}
 \caption{Results for the matrix NSPSDP problem with $X=J_nJ_m^\top$ and $B=H$. Best solution highlighted in bold, if applicable. }
 \label{tab:JH}
 \end{table}
 
 We observe similar results as the ones obtained for the synthetic data. In particular all the methods have the same relative error, ANFGM provides a large norm solution, and FGM is by far the slowest method. In contrast MINGD and CARD perform well and their solutions have low-rank properties. 
 This example shows that it is possible to solve the low-rank NSPSDP problem of large dimension relatively fast, we further discuss the scalability in the next section. 
 
 \reviseLAA{\subsection{A real-world application}
  As a practical example, we study the local compliance estimation problem already considered in \cite{krislock2004local,baghel2022non}, concerning a tiger plush toy. The dimensions of the problem are $n=3$ and $m=12$ and the data matrices (with two digits of accuracy) are
  \[
    X^\top = \begin{pmatrix}
    -0.32 & 0.03 & 0.06 \\
    -0.33 & -0.02 & 0.06 \\
    -0.36 & 0.08 & 0.06 \\
    -0.30 & 0.03 & 0.05 \\
    -0.32 & 0.00 & 0.07 \\
    -0.34 & 0.07 & 0.05 \\
    -0.24 & 0.07 & 0.05 \\
    -0.21 & -0.01 & 0.02 \\
    -0.33 & 0.16 & 0.10 \\
    -0.25 & 0.09 & 0.06 \\
    -0.22 & 0.00 & 0.03 \\
    -0.31 & 0.15 & 0.09
    \end{pmatrix}, \qquad
    B^\top = \begin{pmatrix}
    -1.43 & 0.15 & -0.44 \\
    -1.40 & -0.31 & -0.42 \\
    -1.38 & 0.44 & -0.42 \\
    -1.43 & 0.14 & -0.44 \\
    -1.40 & -0.31 & -0.42 \\
    -1.37 & 0.43 & -0.42 \\
    -1.43 & 0.16 & -0.43 \\
    -1.40 & -0.32 & -0.42 \\
    -1.38 & 0.42 & -0.43 \\
    -1.43 & 0.15 & -0.44 \\
    -1.40 & -0.33 & -0.42 \\
    -1.37 & 0.42 & -0.44
    \end{pmatrix}.
  \]
  The matrix $X$ has rank $r=n=3$, hence the solution to the problem is unique and all the methods converge to the same matrix provided by \cite{baghel2022non,krislock2004local}, with a relative error of $18.99\%$. In order to run the algorithms in a non-full rank setting, we replace $X$ by its best rank-2 approximation, $X_2$. Note that $\|X-X_2\|_F/\|X\|_F = 3.5\%$ so $X_2$ is a good approximation of $X$. Table~\ref{tab:plush} shows the results of the methods applied on $X_2$.
  \begin{table}[ht!]
    \centering
    \begin{tabular}{|c|c|c|c|c|c|}
    \hline
     & Relative error & $\|A\|_F$ & Time (sec.) & $\textnormal{rk}(\symm(A))$ & $\textnormal{rk}(\skeww(A))$ \\
    \hline
    ANFGM (lr) & 0.2022 
    & 6.6785  & \textbf{<0.01} & 2 & 2 \\
    \hline
    FGM (lr)   & 0.2022 
    & 7.4655  & 0.0158 & {\color{red}3} & 2 \\
    \hline
    MINGD (lr) & 0.2022 
    & \textbf{6.6781}  & 0.0143 & 2 & 2 \\
    \hline
    CARD (lr)  & 0.2022 
    & \textbf{6.6781}  & \textbf{<0.01} & 2 & 2 \\
    \hline
    \end{tabular}
    \caption{Results for the tiger plush toy example where $X$ is replaced by its best rank-$2$ approximation, $X_2$.}
    \label{tab:plush}
  \end{table}

  All methods find a solution with minimum relative error. However, CARD and MINGD are the only ones to find a minimizer with minimum norm and minimum rank; FGM does not find a solution with minimum rank nor minimum norm, and ANFGM does not find a solution with minimum norm.  
  }

 \subsection{Scalability} \label{sec:scalability}

 In order to study the scalability of the algorithm, we run the methods considered in the previous examples with increasing dimension and we show the trend of the computational time required to solve the problem. We study two examples, where the dimensions are fixed as $n=m$ and $r=n/2$ for $n=2^j$ with $j=2,\dots,12$. In the first example we consider a sample of $10$ random matrices $X,B\in \R^{n\times m}$ with $\rank(X)=r$, as done for the synthetic data, while for the second example we consider $X=J_n J_n^\top$ and $B=H$, with $J_n$ and $H$ defined as in \eqref{eq:JH}.
 
 \begin{figure}[ht!]
  \centering
   \begin{tabular}{cc}
     \includegraphics[width=0.48\textwidth]{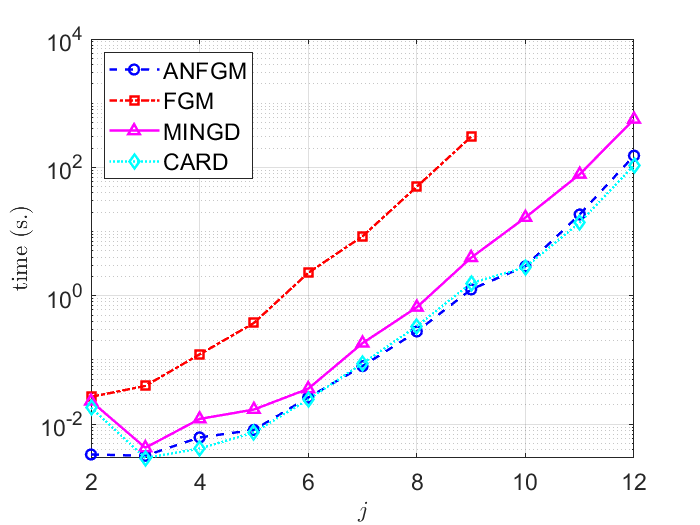}   & \includegraphics[width=0.48\textwidth]{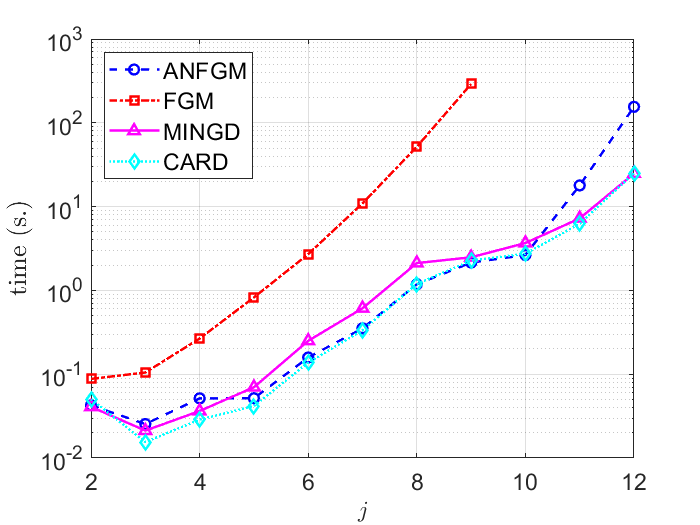} \\
   \end{tabular}
  \caption{Scaling computational times for $n=m=2^j$ and $r=n/2$. On the left the average time over $10$ runs with $X$ a low-rank random matrix, while on the right the example defined by \eqref{eq:JH}. FGM computation exceeded 500 seconds for $n\geq 2^{10}$ in both cases and we did not report its result for larger $n$.}
  \label{fig_scaling}
 \end{figure} 
 
 Figure~\ref{fig_scaling} shows a similar trend for both examples. For ANFGM, MINGD and CARD, the complexity is dominated by the SVD computation ($\io(m n^2)$ flops). The trend of the computational time is approximately a power function of $n$ with exponent between $2$ and $3$. 
 \reviseLAA{Also in these cases ANFGM is generating solutions with very large norms (e.g., larger than $10^{14}$ for $j \geq 9$), as in the previous examples, while FGM is significantly more expensive and its time requirement is larger for $n=2^9$ than that of the other three methods in dimension $n=2^{11}$. 
 
 Interestingly, for the second example (Figure~\ref{fig_scaling}, right plot),  MINGD and CARD are faster than ANFGM for large values of $n$ (namely $n=2^{11}, 2^{12}$). The reason is the initialization step of ANFGM that does not scale as well. 
 For $n = 2^{12}$, the solution of MINGD and CARD coincide, and their CPU time is the same.  
}

 

 \section{Conclusion}
 \label{sec:concl}
 
 In this paper, we proposed a state-of-the-art semi-analytical approach for the NSPSDP problem. Our approach is inspired by that of Baghel et al.~\cite{baghel2022non}, but we resolved an issue when $X$ is rank deficient.  By doing so, we were able to prove that the solution of the NSPSDP problem is always attained (Theorem~\ref{thm:solNSPSDP}). Moreover, we proposed a way to compute the minimum-norm and minimum-rank solution  (Theorem~\ref{thm:minnormNSPSDP}). 
  We illustrated the effectiveness of the new proposed algorithm on several numerical experiments.

 \reviseLAA{
 \section*{Acknowledgment} 
We are very grateful to the editor and the anonymous reviewers who carefully read the manuscript; their feedback allowed us to improve it significantly. 
}
 
 \nocite{*}
 \bibliographystyle{spmpsci}
 \bibliography{biblio}

\end{document}